\documentclass[11pt,twoside]{article}
\usepackage{amsthm}
\usepackage{amssymb,amsbsy,amsmath,amsfonts,amssymb,amscd}
\usepackage{mathrsfs}
\usepackage{graphicx}%
\usepackage{fancyhdr}
\usepackage{color}
\usepackage{algpseudocode}
\usepackage{titlesec}

\theoremstyle{plain} 
\newtheorem{theorem}{Theorem}[section]

\theoremstyle{definition}

\newtheorem{remark}[theorem]{Remark}

\newcommand{\bit}{\begin{itemize}}
\newcommand{\eit}{\end{itemize}}

\newcommand{\ben}{\begin{enumerate}}
\newcommand{\een}{\end{enumerate}}

\newcommand{\ds}{\displaystyle}
\newcommand{\ud}{\, \mathrm{d}}
\newcommand {\eps}  {\varepsilon}
\newcommand {\Neff}  {N_{\text{eff}}}
\newcommand {\Dt}  {\Delta t}

\newcommand {\cell}  {\mathcal C}

\newcommand{\bfx}{{\bf{x}}}
\newcommand{\bfv}{{\bf{v}}}
\newcommand{\bfw}{{\bf{w}}}
\newcommand{\bfn}{{\bf{n}}}
\newcommand{\bfk}{{\bf{k}}}
\newcommand{\bfu}{{\bf{u}}}
\newcommand{\bfE}{{\bf{E}}}
\newcommand{\bfj}{{\bf{j}}}
\newcommand{\R}{\mathbb{R}}
\newcommand{\Sph}{\mathbb{S}}

\newcommand{\be}{\begin{equation}}
\newcommand{\ee}{\end{equation}}
\newcommand{\ba}{\begin{aligned}}
\newcommand{\ea}{\end{aligned}}

\newcommand{\bigo}[1]{\mathcal{O}\left(#1\right)}

\title{A Hybrid Method with Deviational Particles \\ for Spatial Inhomogeneous Plasma\thanks{This research was supported by DOE DE-FG02-13ER26152/DE-SC0010613.} }

\author{Bokai Yan\thanks{Mathematics Department, University of California at Los Angeles,
Los Angeles, CA 90095-1555 USA. byan@math.ucla.edu.}
}

\date{\today}
\begin{document}
\maketitle

\begin{abstract}
In this work we propose a Hybrid method with Deviational Particles (HDP) for a plasma modeled by the inhomogeneous Vlasov-Poisson-Landau system. We split the distribution into a Maxwellian part evolved by a grid based fluid solver and a deviation part simulated by numerical particles. These particles, named deviational particles, could be both positive and negative. We combine the Monte Carlo method proposed in \cite{YC15}, a Particle in Cell method and a Macro-Micro decomposition method \cite{BLM08} to design an efficient hybrid method. Furthermore, coarse particles are employed to accelerate the simulation. A particle resampling technique on both deviational particles and coarse particles is also investigated and improved. This method is applicable all all regimes and significantly more efficient compared to a PIC-DSMC method near the fluid regime.
\end{abstract}

\section{Introduction}
\label{sec:Introduction}

The evolution of a class of collisional plasmas can be modeled by the spatially inhomogeneous Vlasov-Poisson-Landau (VPL) system, with the electromagnetic fields and applied magnetic fields absent,
\begin{subequations}
\label{eq:VPL}
\begin{align}
\label{eq:VPL_f}
&\partial_t f + \bfv\cdot\nabla_\bfx f - \bfE\cdot\nabla_\bfv f= Q_L(f,f), \\
\label{eq:VPL_E}
&-\nabla_\bfx\cdot \bfE = \rho(t,\bfx) = \int f(t,\bfx,\bfv)\ud \bfv.
\end{align}
\end{subequations}
Here $f(t,\bfx,\bfv)$ is the time dependent density distribution function of the charged particles (electrons and ions) in the phase space $(\bfx,\bfv)$, with $\bfx\in\R^3$ the spatial coordinates and $\bfv\in\R^3$ the velocity coordinates.
The left hand side of (\ref{eq:VPL_f}) describes the advection of particles under the electric field $\bfE$, which is solved by the Poisson equation (\ref{eq:VPL_E}). On the right hand side, the Landau (or Landau-Fokker-Planck) operator,
\begin{equation}
  \label{eq:Landau}
  \ba
\ds Q_L(f,f) = \frac{A}{4}\frac{\partial}{\partial v_i} \int_{\R^3} u\sigma_{tr}(u)(u^2\delta_{ij}-u_iu_j)\left(\frac{\partial}{\partial v_j} - \frac{\partial}{\partial w_j}\right) f(\bfw)f(\bfv)\ud \bfw
  \ea
\end{equation}
models the binary collisions due to the long range Coulomb interaction. Here  $u = |\bfv -\bfw|$, $A= 2\pi \left(\frac{e^2}{2\pi\eps_0 m}\right)^2 \log \Lambda$, with $e$ the charge of an individual particle, $m$ its mass, $\eps_0$ the permittivity of free space, and $\log \Lambda$ the Coulomb logarithm.

The small angle collisions dominate in Coulomb interaction. The Landau operator (\ref{eq:Landau}) is derived as the grazing limit of Boltzmann operator (see \cite{book_Cerci} for example) by Landau in 1936 \cite{Landau36}. We refer the reader to Villani \cite{Villani02} and the references therein for a review of the mathematical derivation.  The Landau operator inherits many properties of the Boltzmann operator. First, the density, momentum and energy are conserved during the collision process, i.e.
\[
\left< \phi Q\right> = 0, \quad \text{ for } \phi(\bfv) = 1, \bfv, |\bfv|^2/2,
\]
where $\left< \cdot\right>$ is defined as $\left< \psi\right> = \int_{\R^3} \psi(\bfv) \ud \bfv$.

Furthermore, the entropy is dissipated
\[
\left<  Q \log f \right> \le 0,
\]
with the equal sign holds if and only if $f$ is the local Maxwellian
\be
\label{def:M_local}
M(t,\bfx,\bfv) = \frac{\rho(t,\bfx)}{(2\pi T(t,\bfx))^{3/2}}e^{-\frac{|\bfv-\bfu(t,\bfx)|^2}{2T(t,\bfx)}}.
\ee
Here the density $\rho$, macroscopic velocity $\bfu$ and temperature $T$ are determined by $f$ via
\be
\label{eq:fMsamemoments}
\left< \phi M(t,\bfx,\bfv) \right>  = \left< \phi f(t,\bfx,\bfv) \right>, \quad \text{ for every } \bfx, t,
\ee
with $\phi = 1, \bfv, |\bfv|^2/2$.

The numerical schemes for the VPL system (\ref{eq:VPL}) have been studied extensively. Depending on how to represent the velocity space, these methods can be categorized into two classes, the deterministic (also known as grids based or continuum) methods and probabilistic (also known as particles based) methods. Various techniques have been developed for deterministic methods.
We refer the user to \cite{DLPY15} and the references therein for a recent review.
 However all the deterministic methods are suffered from the {\it{curse of dimensionality}}. The full system (\ref{eq:VPL}) describes a seven dimensional problem (3D in space, 3D in velocity and 1D in time), which leads to tremendous cost of grids based methods.

In contrast, the probabilistic methods, which use particles to represent the velocity space, have the advantage of dimension independence. These methods were first studied by Bird \cite{BirdBook1998}, known as Direct Simulation Monte Carlo (DSMC) method, for rarefied gas. In plasma simulation, Particle in Cell (PIC) methods \cite{Dawson83,book_HE,book_PIC} are widely used in solving the collisionless Vlasov-Poisson system (i.e. system (\ref{eq:VPL}) with $Q_L=0$ in the right hand side). The Landau operator (\ref{eq:Landau}) is solved by Monte Carlo methods, such as Takizuka-Abe method \cite{TakizukaAbe} (henceforth TA) and Nanbu's method\cite{Nanbu97a}. Bobylev and Nanbu \cite{BobylevNanbu}, Bobylev and Potapenko \cite{BP13} proposed a general framework in designing particle methods for (\ref{eq:Landau}). The PIC-DSMC methods, which combines the PIC methods for the advection part and DSMC methods for the collision part, have been widely used in solving the VPL system (\ref{eq:VPL}).

However the particle methods become very inefficient near the fluid regime, where the distribution $f$ is close to the local equilibrium $M$. The reason is that the computation is mainly spent on the collisions between particles sampled from $M$. These collisions have no net effect since $Q_L(M,M)=0$. A hybrid method is more favorable in this regime, which decomposes the distribution as
\be
\label{eq:f_decomp}
f(t,\bfx, \bfv) = M(t,\bfx,\bfv) + f_d(t,\bfx,\bfv).
\ee
The equilibrium $M$ is evolved by a fluid solver and the deviation $f_d$ is solved by particles. Here we allow $f_d<0$ and use positive and negative {\it{deviational}} particles to represent $f_d$.

There have been several related works. In gyrokinetic simulations, the so-called $\delta f$ methods \cite{Byers70, Lee87, DL93, DK95} have been proven to be efficient in collisionless regime. Hybrid methods with weighted deviational particles have also been designed when the collisions are modeled by BGK operator \cite{CCL12} and linear Fokker-Planck operator \cite{SWHK15}. However all these methods are very difficult to efficiently handle binary collisions. Recently
Caflisch et al. \cite{CWDCD08} and Ricketson et al. \cite{RRCD14} studied the thermalization/dethermalization methods with $f_d(\bfv)\ge0$ for the spatial homogeneous Coulomb collisions. However they are difficult to be generalized to the inhomogeneous system (\ref{eq:VPL}) without allowing $f_d$ to be negative.
For rarefied gas, Hadjiconstantinou and co-workers \cite{BakerHadjiconstantinou05, Wagner08, BakerHadjiconstantinou08, HomolleHadjiconstantinou_PhysFluid, HomolleHadjiconstantinou_JCP, PLH14} have studied a low-variance deviational simulation Monte Carlo (LVDSMC) method in which $f_d$ could be negative. However, as discussed in Section \ref{sec:collision_Landau}, the extension to long range Coulomb collisions leads to dramatic growth of particle numbers and hence a numerical instability.


A different Monte Carlo method with deviational particles for Coulomb collisions has been designed by the authors in \cite{YC15}. By introducing coarse particles as a rough approximation of $f$, a large number of collisions are combined efficiently and the growth of particle numbers is reduced from the {\it{numerical}} scale to {\it{physical}} scale. In this work we seek to extend the method in \cite{YC15} to the inhomogeneous VPL system (\ref{eq:VPL}).

The goal of this work is to design a hybrid method with deviational particles which (a) is applicable in all regimes and (b) is much more efficient than the PIC-DSMC method near the fluid regime. The starting point is an operator splitting between the collision part and the advection part. The collision part is solved as in the homogeneous case cell by cell. As for the advection part, we combine a macro-micro projection method \cite{BLM08} and a Particle-In-Cell method to evolve the macroscopic quantities and all particles. As a result, the number of particles grows in the physical scale and a particle resampling method is applied on both deviational and coarse particles to further control the growth.

The rest of the article is organized as follows. In Section \ref{sec:method_whole} we develop hybrid methods with deviational particles for the VPL system, as well as a simplified system, the VP-BGK system. We describe the algorithms in solving the collision part in Section \ref{sec:collision} and the advection part in Section \ref{sec:transport}. Next in Section \ref{sec:resample} we investigate and improve the particle resampling technique to control the particle growth. Finally in Section \ref{sec:numerics} various numerical simulations are performed to demonstrate the accuracy and efficiency of our hybrid methods.

\section{Hybrid Methods with Deviational Particles}
\label{sec:method_whole}

As in (\ref{eq:f_decomp}), we decompose the distribution $f(t,\bfx, \bfv)$ into a Maxwellian part $M(t,\bfx,\bfv)$ and a deviation $f_d(t,\bfx,\bfv)$, which can be both positive and negative.


We use a hybrid representation of $f$. A uniform grid is placed in each spatial direction. The Maxwellian part $M(t,\bfx,\bfv)$ is given by (\ref{def:M_local}), with the macroscopic quantities defined on the grids.

For simplicity, the Maxwellian part is chosen to have the same the moments of $f$ as in (\ref{eq:fMsamemoments}). These moments are solved by a fluid solver.

The deviation part is represented by particles. As $f_d$ could be both positive and negative, we use {\it{deviational particles}} in the representation. In other words, in each spatial cell we sample and evolve {\it{positive}} particles for the positive part $(f_d)_+$ and {\it{negative}} particles for the negative part $(f_d)_-$. We denote $N_{p,k}$, $N_{n,k}$ and $N_{d,k}$ the number of positive particles, negative particles and all deviational particles in the $k$th spatial cell $\cell_k$. Then $N_{d,k} = N_{p,k} + N_{n,k}$, for every $k$.

The total moments of the deviational particles in each spatial cell are 0, i.e., for $\phi = 1, \bfv, |\bfv|^2/2$,
\be
\label{eq:fd_zeromoments}
\left<\phi f_d(x_k)\right> = \Neff \sum_{\bfv_p\in\cell_k} \phi(\bfv_p) - \Neff \sum_{\bfv_n\in\cell_k} \phi(\bfv_n) = 0 , \quad \text{ for every } k \text{ and } t.
\ee
Here in the summation $\bfv_p$ are over all positive particles and $\bfv_n$ are over all negative particles in each spatial cell $\cell_k$. $\Neff$ is the effective number which represents the physical densities per numerical particle. In this work we choose $\Neff$ to be a prescribed constant for all the time. In other words, we work on the equally weighted particles.

In practice, (\ref{eq:fd_zeromoments}) is only valid in expectation. It cannot hold exactly due to statistical fluctuation unless one enforces the conservations after each step.


To solve the VPL system (\ref{eq:VPL}), we start with an operator splitting between the collision part and the advection part. In the following we describe how to evolve the two components $M$ and $f_d$ in the two substeps.

\subsection{Deviational Particles in Collisions}
\label{sec:collision}

In the substep of collision, the Maxwellian part $M$ is invariant. Here we describe how to evolve the deviation part $f_d$ for different collisions.

\subsubsection{BGK operator}

Let us first consider a simplified model called the Vlasov-Poisson-BGK (VP-BGK) system, which is the same as the VPL system (\ref{eq:VPL}) except the collision term $Q_L$ is replaced by the BGK operator,
\be
\label{eq:BGK}
Q_{BGK}(f) = \mu (M-f).
\ee
This operator models the relaxation of the distribution $f$ to its local equilibrium $M$. Here the collision frequency $\mu$ is independent of $\bfv$. Inserting the decomposition (\ref{eq:f_decomp}) in the collision substep $\partial_t f = Q_{BGK}(f)$, one has
\[\partial_t f_d = -\mu f_d.\]
For this linear model, there are no interactions between particles.

A Monte Carlo method can be designed by applying a forward Euler method,
 \be
 \label{eq:BGK_exp}
  f_d^{n+1} = (1 - \mu\Dt) f_d^n.
  \ee
Hence each (positive and negative) particle is removed with probability $\mu\Dt$ in one time step. Clearly one needs $\mu\Dt\le1$.

One can also apply a backward Euler method,
 \be
  \label{eq:BGK_imp}
  f_d^{n+1} = \frac{1}{1+\mu\Dt} f_d^n.
  \ee
Again a Monte Carlo method can be designed accordingly. Each (positive and negative) particle is removed with probability $\frac{\mu\Dt}{1+\mu\Dt}$ in one time step. In this case there is no restriction on the time step $\Dt$.

\subsubsection{Landau operator}
\label{sec:collision_Landau}

Now we consider the Coulomb collisions modeled by the bilinear Landau operator (\ref{eq:Landau}). The binary interaction implies that one has to include all six types of collisions: P-P, P-N, N-N, P-M, N-M and M-M. Here, for example, a P-N collision means a collision between a {\it{positive}} particle and a {\it{negative}} particle. M means the Maxwellian part. Note that one can omit the M-M collisions which have no net effect. This saves the major cost compared to a full particle method.

To simulate the other five types of collisions, one has two methods: (a) simulate each type of collisions individually; (b) group different types of collisions and simulate them together.

\

The first method is straightforward. One only needs to determine the numbers of different types of collisions and formulate corresponding collision rules. These rules have been derived in simulating the binary collisions in {\it{rarefied gas}} which consists of charge neutral particles (see \cite{BakerHadjiconstantinou05}). They can be summarized as
 \begin{equation}
  \label{eq:boltzmann_rule}
  \begin{aligned}
    \ds \mbox{P-P:}& \quad \bfv_+, \bfw_+ \to \bfv_+', \bfw_+', \\
    \ds \mbox{P-N:}& \quad \bfv_+, \bfw_- \to 2\bfv_+, \bfv_-', \bfw_-', \\
    \ds \mbox{N-N:}& \quad \bfv_-, \bfw_- \to 2\bfv_-, 2\bfw_-, \bfv_+', \bfw_+', \\
    \ds \mbox{P-M:}& \quad M, \bfv_+ \to M, \bfw_-, \bfv_+', \bfw_+', \\
    \ds \mbox{N-M:}& \quad M, \bfv_- \to M, \bfw_+, \bfv_-', \bfw_-'.
  \end{aligned}
\end{equation}
Here we explain the P-N collision in detail. The other rules can be explained similarly. In a P-N collision between a positive particle $\bfv_+$ and a negative particle $\bfw_-$, the two incident particles are removed, four new particles (two positive particles with velocity $\bfv$ and two negative particles with velocities $\bfv'$ and $\bfw'$) are created. Here the post-collisional velocities $\bfv'$ and $\bfw'$ are defined in the usual way in Boltzmann equation (see\cite{book_Cerci} for example). We refer the reader to \cite{YC15} for more details on the derivation and implementation of this method.

One can apply these collision rules to the Coulomb gas which consists of charged particles. However a problem of the method based on (\ref{eq:boltzmann_rule}) is that the particle number grows due to collisions in each time step. The increment depends on the number of collisions in one step. For rarefied gas, this problem may be manageable since the number of collisions in one time step is $\bigo{\Dt}$, hence $N_d^{t+\Dt} = (1+c\Dt) N_d^t$. As $\Dt\to0$, the number of particles is $N_d^t = N_d^0 e^{ct}$, which is $\Dt$ independent. This increase is in the {\it{physical}} time scale.

However the extension of this method to Coulomb collisions suffers a more severe increase due to the fact that in Coulomb gas all particles are involved in collisions no matter how small the time step $\Dt$ is. This is the effect of long range interaction between charged particles. As a result, the particle number increases in the {\it{numerical}} time scale, i.e. $N_d^{t+\Dt} = c N_d^t$, with $c\ge1$ independent of $\Dt$. Hence the total number of particles increases in a {\it{numerical}} scale.

\

Therefore, instead of simulating each collisions individually, we seek to group different types of collisions and simulate them together to reduce the number of collisions in Coulomb gas.
More precisely, in \cite{YC15} we proposed to split the binary collision operator as
\[ Q(f, f) = Q(f,f_d) + Q(f, M) = Q(f,f_d) + Q(f_d, M),\]
and the spatial homogeneous equation $\partial_t f = Q(f,f)$ becomes
\begin{equation}
\label{eq:systemtosovle_coll_nocoarse}
  \left\{\ba
  \partial_t M & =0, \\
  \partial_t f_d & = Q(f, f_d) + Q(f_d, M). \\
  \ea\right.
\end{equation}
Here the operator $Q(f,f_d)$ is defined to be asymmetric in $f$ and $f_d$, so that $Q(f, f_d)$ represents the changes in $f_d$ due to collisions with $f$.
Note that $Q(f, f_d) = Q(f, (f_d)_+) - Q(f, (f_d)_-)$. One has grouped M-P, P-P and N-P collisions in the first term $Q(f, (f_d)_+)$ since
 $$Q(f, (f_d)_+) = Q(M, (f_d)_+) + Q((f_d)_+, (f_d)_+) - Q((f_d)_-, (f_d)_+).$$
 Similarly M-N, P-N and N-N collisions have been grouped and represented by the term $Q(f, (f_d)_-)$.

To solve the term $Q(f,f_d)$, one performs regular collisions between particles sampled from $f$ and the deviational particles. However, sampling from $f$ at every time step is very expensive, since one needs to recover the distribution $(f-M)$ in 3D velocity space from the deviational particles. Instead of sampling from $f$, we proposed in \cite{YC15} to include a bunch of {\it{coarse particles}} which solve the original homogeneous equation $\partial_t f = Q(f,f)$. To sample a particle from $f$, one only needs to randomly pick up a coarse particle. More specifically, we proposed to solve the following system
\begin{equation}
\label{eq:systemtosovle_coll}
  \left\{\ba
  \partial_t f_c & = Q(f_c, f_c), \\
  \partial_t M & =0, \\
  \partial_t f_d & = Q(f_c, f_d) + Q(f_d, M). \\
  \ea\right.
\end{equation}
Here $f_c$ solves the original homogeneous Landau equation. However, $f_c$ serves as a coarse approximation of the solution and is represented by a relatively small number of particles, i.e. the {\it{coarse particles}}. One can evolve the coarse particles by a regular Monte Carlo method, for example the TA method or Nanbu's method.

Now the evolution of the deviational particles involve two steps. First, to solve the $Q(f_c, f_d)$ part, the deviational particles are updated by performing regular collisions with the coarse particles. In this step $N_d$, the number of deviational particles, does not change. Then, a number of $\bigo{\Delta t N_d}$ deviational particles are sampled from the distribution $Q(f_d, M)$ and added into the system. Therefore after one step of collision, $N_d$ increases by $N_d^{t+\Dt} = (1+c\Dt) N_d^t$, i.e., in the {\it{physical}} time scale.

The Maxwellian part $M$ does not change in the collision step. Therefore the coarse particles provide a coarse, direct simulation of $f$, while the deviational particles give a fine, deviational simulation of $(f-M)$. This method can be applied to any type of binary collisions, for example, the short range hard sphere collision and the long range Coulomb collision.

To perform the collisions between the deviational particles and the coarse particles, i.e. the $Q(f_c, f_d)$ part, the number of coarse particles should be large enough. One needs
\be
\label{eq:NfNpNn}
N_{c,k}\ge N_{d,k}, \quad \text{ for every } k,
\ee
with $N_{c,k}$ the number of coarse particles in the spatial cell $\cell_k$. We assume all the coarse particles in all cells carry the same weight to avoid the collisions between particles with variable weights. However, the effective number of the coarse particles could be time dependent as $N_d$ varies in time. We will discuss more on this in Section \ref{sec:resample_coarse}.

Finally,  note that each new sampled particle needs both spatial and velocity coordinates. To sample a particle with velocity distribution to be $Q(f_d, M)$, we refer the reader to Section 4 in \cite{YC15}. The spatial coordinate of each sample is set to be uniformly distributed in current spatial cell $\cell_k$.

\begin{remark}
  The system (\ref{eq:systemtosovle_coll}) describes an particle method for a general bilinear operator $Q(f,f)$ which admits an equilibrium $M$, i.e., $Q(M,M) = 0$. One can formally apply it to Landau operator (\ref{eq:Landau}). However, the existing Monte Carlo methods for Coulomb collisions such as Takizuka-Abe method \cite{TakizukaAbe} and Nanbu's method\cite{Nanbu97a} are not solving the Landau equation $\partial_t f = Q_L(f,f)$ directly. They instead solve the formula proposed in \cite{BobylevNanbu},
  \begin{equation}
  \label{eq:LFP_approxi}
  \ds f(\bfv,t+\Delta t)=\frac{1}{\rho} \iint_{\R^3\times \Sph^2} D\left(\frac{\bfu\cdot\bfn}{u}, A \frac{\Delta t}{u^3}\right) f(\bfw',t) f(\bfv',t) \ud \bfw\ud \bfn,
\end{equation}
with $D(\cdot,\cdot)$ an approximated kernel determined by the numerical algorithms, $\bfu = \bfv - \bfw$, $\rho$ the density of $f$ and $A$ the coefficient in (\ref{eq:Landau}). This is a first order approximation of the solution $f(t^{n+1})$.
Therefore, in practice one applies the method  (\ref{eq:systemtosovle_coll}) to (\ref{eq:LFP_approxi}), leading to a first order approximation of (\ref{eq:systemtosovle_coll}),
\be
\label{eq:systemtosovle_Landau}
\left\{
\ba
\ds f_c(\bfv,t+\Delta t) &=  \frac{1}{\rho} \iint D \, { f_c(\bfw')} \, {f_c(\bfv')}  \ud \bfw\ud \bfn, \\
\ds M(\bfv,t+\Delta t) &= M(\bfv,t), \\
\ds f_d(\bfv,t+\Delta t) &=  \frac{1}{\rho} \iint D \, { f_c(\bfw')}\, {f_d(\bfv')}  \ud \bfw\ud \bfn + \Delta M(\bfv),
\ea
\right.
\ee
where
\be
\label{eq:distri_sample}
\ds \Delta M(\bfv) = \frac{1}{\rho} \iint D \,f_d(\bfw') M(\bfv')\ud \bfw\ud \bfn- \frac{\rho_p - \rho_n}{\rho}M(\bfv)
\ee
We refer the reader to Section 4 in \cite{YC15} for details.
\end{remark}

\subsection{Deviational Particles in Advection}
\label{sec:transport}

Now consider the advection part, which is commonly known as the collisionless Vlasov-Poisson equation,
\be
\label{eq:VP}
\partial_t f + \bfv\cdot\nabla_\bfx f - \bfE\cdot \nabla_{\bfv} f = 0,
\ee
with the electric field $\bfE$ given by the Poisson equation (\ref{eq:VPL_E}), where the density $\rho$ is computed from
\[
\rho(x_k) = \rho_M(x_k) + \frac{\Neff}{|\cell_k|} (N_{p,k} - N_{n,k})
\]
in each spatial cell $\cell_k$. Here $\rho_M$ is the density of the Maxwellian part, $N_{p,k}$ and $N_{n,k}$ are the numbers of positive and negative deviational particles in that cell. $|\cell_k|$ is the size of the cell. Note that $N_{p,k}$ and $N_{n,k}$ might not be exactly the same due to the statistical fluctuation.

In our method, $\rho_M$, $N_p$ and $N_n$ are computed explicitly from $M$ and deviational particles at the beginning of the advection step. Then $\bfE$ is computed from the Poisson equation (\ref{eq:VPL_E}). Denote
\[\mathcal T = \bfv\cdot\nabla_{\bfx} - \bfE\cdot \nabla_\bfv.\]
Now $\mathcal T$ is a linear operator during current time step.

In the advection step, we need to evolve both the Maxwellian part $M$ and the deviational particles which represent $f_d$. For the Coulomb collisions, we also need to evolve the coarse particles which represent $f_c$.

For $M$ and $f_d$, we employ the macro-micro decomposition introduced in \cite{BLM08}. First (\ref{eq:VP}) is reformulated as
\be
\label{eq:Boltzmann_mm_MPN}
\left\{
\ba
& \frac{\partial}{\partial t} \left< M \phi \right> + \nabla_\bfx \cdot \left< \bfv M \phi\right> + \nabla_\bfx \cdot \left< \bfv f_d \phi\right> = \left<
\bfE\cdot \nabla_{\bfv} f\phi\right>  = (0, -\rho \bfE, -\rho \bfu \cdot \bfE)^T, \\
& \frac{\partial f_d}{\partial t} + \mathcal T f_d = -\left(\frac{\partial}{\partial t} + \mathcal T \right) M, \\
\ea
\right.
\ee
with $\phi = 1, \bfv, |\bfv|^2/2$. Here the first equation is obtained by taking the moments of both sides of (\ref{eq:VP}), while the second equation is just a rearrangement of (\ref{eq:VP}).

To solve the first equation, one can apply a fluid solver to compute the term $\nabla_x \cdot \left< \bfv M \phi\right>$. For the next term one can use the approximation
\be
\label{eq:vfphi}
\left< \bfv f_d \phi\right> |_{x_k} \approx \Neff \sum_{\bfv_p\in\cell_k} (\bfv_p \phi(\bfv_p)) - \Neff \sum_{\bfv_n\in\cell_k} (\bfv_n \phi(\bfv_n)).
\ee

For the second equation of (\ref{eq:Boltzmann_mm_MPN}), to eliminate the time derivative term $\partial_t M$, we define a projection operator $\Pi_M$ by
\be
\label{eq:projection}
\Pi_M \psi = \frac{M}{\rho_M} \left[\left<\psi\right> + \frac{(\bfv-\bfu_M)\cdot\left<(\bfv-\bfu_M)\psi\right>}{T_M} + \frac{1}{2d}\left(\frac{|\bfv-\bfu_M|^2}{T_M} - d\right)\left<\left(\frac{|\bfv-\bfu_M|^2}{T_M} - d\right)\psi\right>\right],
\ee
where $\rho_M$, $\bfu_M$ and $T_M$ are the density, macroscopic velocity and temperature  of $M$. $d$ is the dimension of velocity space. We take $d = 3$ in this work.

From (\ref{eq:fd_zeromoments}), one has
\[ \Pi_M f_d = \Pi_M (\partial_t f_d) = (I-\Pi_M)(\partial_t M) =0.\]
Apply the operator $(I-\Pi_M)$ on both sides of the second equation of (\ref{eq:Boltzmann_mm_MPN}), one obtains
\[
\frac{\partial f_d}{\partial t} + \mathcal T f_d = -(I-\Pi_M)\left(\mathcal T M \right) + \Pi_M (\mathcal T f_d).
\]

This equation governs the evolution of the deviational particles.
There are two processes involved: the advection of particles and the sampling from the source term on the right hand side. The advection process is just like a regular PIC method: one moves all the deviational particles by solving the ODE system
 \be
 \label{eq:VP_ODE}
  \dot \bfx = \bfv, \quad \dot \bfv = -\bfE,
 \ee
 with $\bfE$  the electric field in the spatial cell the particle resides. In this process $N_d$ is not changed. Then in the sampling process, deviational particles are created due to the fact that $M$ is driven away from the shape of a Maxwellian after convected. Noting that $f_d$ is not included in the first term, a number of $\bigo{\Dt \left(\Neff^{-1} + N_d\right)}$ particles are created. The deviational particles are created even if they do not exist at beginning. After created, $N_d$ increases in the {\it{physical}} time scale, just as in the Coulomb collision step (\ref{eq:systemtosovle_coll}).  In Appendix \ref{append:advect_source} we give an efficient algorithm to sample velocities from the source term.

Again the spatial coordinate of each new sampled particle is set to be uniformly distributed in current spatial cell $\cell_k$.

Furthermore, note that the moments of the deviational particles remain zero in expectation. This is due to the fact that, for any distribution $\psi$,
\[ \left<\phi (I-\Pi_M) \psi\right>=0, \quad \text{ for } \phi = 1, \bfv, |\bfv|^2/2.\]

Finally, for the VPL system, one also needs to transport the coarse particles which represent $f_c$ in (\ref{eq:systemtosovle_coll}). Similar to the collision step,  $f_c$ just solves the original Vlasov-Poisson equation (\ref{eq:VP}) in the advection step. One can simply evolve the coarse particles by a regular PIC method.

Therefore the system to be solved can be summarized as,
\begin{subequations}
\label{eq:systemtosovle_adve}
\begin{align}
\label{eq:advection_fcoarse}
&\frac{\partial f_c}{\partial t} + \mathcal T f_c = 0, \\
\label{eq:advection_M}
& \frac{\partial}{\partial t} \left< M \phi \right> + \nabla_\bfx \cdot \left< \bfv M \phi\right> + \nabla_\bfx \cdot \left< \bfv f_d \phi\right>  = (0, -\rho \bfE, -\rho \bfu\cdot \bfE)^T, \\
\label{eq:advection_fd}
&\frac{\partial f_d}{\partial t} + \mathcal T f_d = -(I-\Pi_M)\left(\mathcal T M \right) + \Pi_M (\mathcal T f_d).
\end{align}
\end{subequations}

\begin{remark}
  The advection of deviational particles was also studied in \cite{HomolleHadjiconstantinou_JCP}. In that work, the Maxwellian part $M$ is assumed to be invariant during the advection time step. The time evolution of the macroscopic quantities are left in the collision part. Hence in the advection step one only updates the deviational particles, by solving the second equation in (\ref{eq:Boltzmann_mm_MPN}) with $\frac{\partial}{\partial t} M=0$. The macro-micro projection was not used.
\end{remark}

\subsection{Summarize}

We refer to the Hybrid methods with Deviational Particles developed in this section as the HDP methods. Now we summarize the HDP methods for VP-BGK system and VPL system.

First, for both systems, an operator splitting is applied between the advection part and the collision part. Then different algorithms are applied in each substep.

\

{\bf{HDP method for VP-BGK system}}.
\begin{itemize}
  \item In the collision substep, one solves (\ref{eq:BGK_exp}) or (\ref{eq:BGK_imp}).
  \begin{itemize}
    \item The Maxwellian part $M$ keeps invariant.
    \item The deviational particles are removed with rates $\mu\Delta t$ or $\frac{\mu\Delta t}{1+\mu\Delta t}$.
  \end{itemize}
  \item In the advection substep, one solves (\ref{eq:advection_M}) and (\ref{eq:advection_fd}).
    \begin{itemize}
    \item The moments of the Maxwellian part $M$ are updated by applying a fluid solver on (\ref{eq:advection_M}).
    \item The deviational particles move in the phase space by solving (\ref{eq:VP_ODE}). New deviational particles sampled from the right hand side of (\ref{eq:advection_fd}) are added into the system.
  \end{itemize}
\end{itemize}

{\bf{HDP method for VPL system}}. Coarse particles are needed to solve this system.
\begin{itemize}
  \item In the collision step, one solves (\ref{eq:systemtosovle_coll}).
  \begin{itemize}
    \item The coarse particles are updated by applying a regular DSMC method to simulate the Coulomb collision.
    \item The Maxwellian part $M$ keeps invariant.
    \item The deviational particles collide with the coarse particles. New deviational particles sampled from $\Delta t  Q(f_d, M)$ are added into the system.
  \end{itemize}
  \item In the advection step, one solves (\ref{eq:systemtosovle_adve}).
  \begin{itemize}
    \item The coarse particles move in the phase space by solving (\ref{eq:VP_ODE}). This is the traditional PIC method.
    \item The moments of the Maxwellian part $M$ are updated by applying a fluid solver on (\ref{eq:advection_M}).
    \item The deviational particles move in the phase space by solving (\ref{eq:VP_ODE}). New deviational particles sampled from the right hand side of (\ref{eq:advection_fd}) are added into the system.
  \end{itemize}
\end{itemize}

\begin{remark}
\label{remark:growth}
In the VP-BGK system, the collision step automatically reduces $N_d$  by a constant ratio in one time step. As a result, for a moderate collision frequency $\mu$, the collision step can stabilize the advection step. Furthermore, the positive and negative deviational particles cannot be mixed in neither the collision nor the advection step.

In the VPL system, the number of deviational particles increases by $\bigo{\Delta t N_d}$ in both advection step and the collision step. Hence $N_d$ grows exponentially in the physical time scale. However, as will be discussed in next section, the Coulomb collision can mix the positive and negative deviational particles, making it possible to reduce $N_d$ by canceling particle pairs if they are close in the phase space.
\end{remark}




\section{Control of Particle Growth}
\label{sec:resample}

As discussed in Remark \ref{remark:growth}, the number of deviational particle grows exponentially in the physical scale when solving the VPL system (\ref{eq:VPL}) by the HDP method. A possible way to control the growth is to remove particle pairs if they are close in the phase space. Similar to the homogeneous case studied in \cite{YC15}, it is both necessary and efficient to apply a particle reduction technique.

\begin{itemize}
  \item First, to perform the collisions between the deviational particles and the coarse particles, i.e. the $Q(f_c,f_d)$ term in (\ref{eq:systemtosovle_coll}), one needs (\ref{eq:NfNpNn}) to be valid in each cell. However, the coarse particles are evolved by applying a regular DSMC method on the collisions and a regular PIC method on the advection. Hence $N_c$ does not change while $N_d$ increases in the physical scale. After a moderate time period, the assumption (\ref{eq:NfNpNn}) could fail. Therefore it is necessary to reduce $N_d$ periodically.
  \item Second, the first term $Q(f,f_d)$ in the collision step drives both $(f_d)_+$ and $(f_d)_-$ to the same distribution. The positive and negative deviational particles are approaching each other. Therefore by canceling particle pairs which are close in phase space, one can reduce the number of particles and expect $N_d$ to be bounded and even reduced. This increases the efficiency of the whole method.
\end{itemize}

\subsection{Resample Deviational Particles}

A particle resampling method was studied in \cite{YC15}. One first finds an approximation of the explicit formula of the deviation part $f_d$ and then resamples deviational particles from it. More specifically, $f_d$ is approximated by
\be
\label{eq:resample_formula}
\begin{aligned}
f_{d}(\bfv) &= \sum_\bfk\left< f_d, \phi_\bfk \right> \phi_\bfk(\bfv)  \approx \sum_{|\bfk|<K}\left< f_d, \phi_\bfk \right> \phi_\bfk(\bfv) \\
  & \approx  \Neff\sum_{|\bfk|<K} \left( \sum_{\bfv_{p}} \phi_\bfk (\bfv_{p}) - \sum_{\bfv_{n}} \phi_\bfk (\bfv_{n}) \right) \phi_\bfk(\bfv).
\end{aligned}
\ee
Here $\{\phi_\bfk\}$ is a complete basis set. $\left<\cdot, \cdot\right>$ is the inner product, and $K$ is a cutoff parameter.

The total cost is $\bigo{N_dK^3}$ to evaluate the coefficients in (\ref{eq:resample_formula}) and to resample new particles. We find at least $K=20$ is needed to make a resampling with reasonable accuracy. This makes the cost to be $\bigo{8000N_d}$ which is pretty large. This cost is acceptable in the spatial homogeneous case, since the distribution decays to equilibrium monotonically and the number of deviational particles becomes very small as time evolves. But in the spatial inhomogeneous case, the advection part also creates a lot of particles and drives the distribution away from the Maxwellian. A lot more resamplings are needed. To accelerate the resampling, in Appendix \ref{append:resample_acc} we propose an approximation of the formula (\ref{eq:resample_formula}), which reduces the total cost to $\bigo{N_d + K^3\log K}$.

In the spatial homogeneous case, the resampling is performed when the condition (\ref{eq:NfNpNn}) is violated. So a straightforward algorithm for the inhomogeneous system is that, one performs the resamping in the cell where (\ref{eq:NfNpNn}) is violated. However this strategy is not very efficient. Consider the case that the condition (\ref{eq:NfNpNn}) is violated and the resampling is performed in some cell for the first time. Then the neighboring spatial cells, in which no resampling has been performed and hence many deviational particles are accumulated, might convect many particle pairs into this cell in the next time step. This makes the previous resampling meaningless.

Here we use the following strategy. Whenever the condition (\ref{eq:NfNpNn}) is violated in one cell, the resampling of deviational particles is performed in all the cells $\cell_k$ where
\be
\label{eq:NfNpNn_beta}
N_{d,k} \ge \beta N_{c,k},
\ee
with a prescribed constant $\beta<1$.

The case $\beta=1$ is exactly the algorithm that the resampling is performed only in the cells violating (\ref{eq:NfNpNn}). If $\beta$ is very small, many unnecessary resamplings will be performed in the cells where $N_{d,k} \ll N_{c,k}$. In practice we find a simple choice $\beta = 0.9$ gives quite satisfactory efficiency.


\begin{remark}
  In the collision process, the two types of particles are mixed due to the collisions with the coarse particles.
  However there is no such effect in the advection step . Once the electric field is computed, the advection operator is linear and time invertible. The supports of the positive and negative particles keep disjoint during the advection, if they are disjoint initially. Therefore to make effective resampling, the collision cannot be weak. This is also the starting point of using deviational particles, i.e. we assume the distribution is close to the local equilibrium.
\end{remark}

\begin{remark}
  For the spatial inhomogeneous case, we need to apply the particle resampling in each spatial cell, instead of the whole $(\bfx,\bfv)$ phase space, due to two reasons. First, the error introduced by the cutoff after $K$ terms in (\ref{eq:resample_formula}) depends on the smoothness of the underlying function $f_d$. In velocity direction $f_d$ is smooth due to the Coulomb collisions (see \cite{Safadi07}). However $f_d$ might be discontinuous in space coordinates due to boundary conditions or external forces.  Second, resampling in the whole $(\bfx,\bfv)$ phase space leads to a memory usage of $\bigo{K^3N_x^3}$ to store the coefficients in (\ref{eq:resample_formula}). This loses the advantages of using particles.
\end{remark}

\subsection{Resample Coarse Particles}
\label{sec:resample_coarse}

The coarse particles also need resampling, which means changing the effective number of the coarse particles. This effective number, denoted by $\Neff^c$, could be both increased and decreased by removing or adding some coarse particles. Noting that all the coarse particles in all cells carry the equal weight, we have to perform resampling in all cells once it is needed in one cell.

Let $N_c$ be the total number of coarse particles. There are three cases one needs to resample the coarse particles.

{\bf{Case 1. Reduce $N_c$ to improve efficiency. }}

We first recall the resampling of coarse particles in the spatial homogeneous case. Note that $(f_d)_+$ and $(f_d)_-$, the positive and negative parts of the deviation, approach the same equilibrium as time evolves. Hence $N_d$, the number of deviational particles, is reduced significantly after each particle resampling. As the only requirement on the coarse particle is (\ref{eq:NfNpNn}), one can resample the coarse particles by simply discarding a fraction of them. In this case the effective number $\Neff^c$ increases. Eventually only a small number of deviational and coarse particles are needed.

However this is not true in the spatial inhomogeneous case. The advection part drives the distribution away from the local equilibrium. The balance between the advection part and the collision part keeps $N_d$ from being diminished. $N_d$ is not likely to be reduced significantly as in the homogeneous case. Since the resampling of the coarse particles always involves all the cells, in this work we choose {\it{not}} to reduce $N_c$ even if $N_d$ is reduced through resampling in some cell.

{\bf{Case 2. Increase $N_c$ if the resampling of deviational particles fails.}}

The resampling of deviational particles may fail in some cell if the growth of $N_d$ is mainly due to the advection part. This could happen, for example, in the early stage of simulation with initial distribution to be local equilibrium, or in the regime with large spatial variation. In these cases, the overlap between the positive and negative deviational particles is very small, leading to very small and even no reduction of $N_d$ in a resampling. To satisfy (\ref{eq:NfNpNn}), we need to introduce more coarse particles.

However, one cannot obtain more {\it{independent}} coarse particles from only the original coarse particles. To increase $N_c$ effectively, one needs to sample directly from $f=M+f_d$, with $f_d$ recovered from the deviational particles as in (\ref{eq:resample_formula}) or in Appendix \ref{append:resample_acc}. One needs to sample more coarse particles in all cells to ensure that their effective numbers remain the same.

{\bf{Case 3. To conform to the fine solution.}}

After simulation of a relatively long time without resampling coarse particles, the accumulated error in the coarse solution $f_c$ might become significant. In this case one can discard the original coarse particles and resample them from $f=M+f_d$ directly. This keeps the coarse solution consistent with the fine solution.

In practice, this resampling includes two steps. First, we resample the deviational particles in all cells. Then coarse particles are resampled with the new effective number
\be
\label{eq:newNeff}
(\Neff^c)_\text{new} = \max \left(\min_k \frac{\rho_M(x_k)|\cell_k|}{\gamma N_{d,k}}, (\Neff^c)_\text{old}\right),
\ee
with $\gamma \gtrsim 1$ a prescribed constant. We take $\gamma = 1.1$ in the simulations. This resamples no more than the old coarse particles without violating the condition (\ref{eq:NfNpNn}) in each cell.

\

As a result, the resampling of coarse particles in case 2 increases $N_c$. However,  as time evolves the distribution $f$ might approach some global equilibrium, which demands less deviational particles. The resampling in case 3 has a chance to reduce the number of coarse particles.

%

\subsection{Assign x coordinates after resampling}

Finally, one needs to assign an x coordinate to each resampled particle. A simple choice is to pick up a random number uniformly distributed in the spatial cell as the new coordinate. This introduces a first order error in cell size $\Delta x$. It is possible to obtain a second order reconstruction by approximating the coefficients in the expansion (\ref{eq:resample_formula}) by linear functions in $x$. This will be further investigated in future work.

\subsection{Another viewpoint of the coarse particles}

The evolution and resampling of coarse particles can be explained in another way. As discussed in Section \ref{sec:collision_Landau}, sampling from $f$ at every time step to solve (\ref{eq:systemtosovle_coll_nocoarse}) is very expensive. Hence once particles are sampled from $f = M + f_d$, they are not discarded immediately after that time step. Instead these particles are evolved and reused in the following time steps. After a while, when these particles cannot represent $f$ accurately due to the error accumulation, we discard these particles and sample a new collection of particles from $f = M + f_d$. This is exactly the evolution and resampling of coarse particles.

\section{Numerical Tests}
\label{sec:numerics}

In the following we perform numerical simulations on the VP-BGK system and the VPL system (\ref{eq:VPL}).
The HDP methods are applicable to phase spaces with any dimensions. Here for simplicity our tests are performed on the 1D spatial and 3D velocity coordinates.

In the tests, we take spatial domain to be $x\in[0,4\pi]$, with a periodic boundary condition in $x$. Note that a truncation in velocity space is not necessary since we use particles to represent the distribution in $\bfv$. In our tests, we take the coefficient $A = 10$ in (\ref{eq:Landau}). Unless specified, we always use $N_x=400$ grid points in $x$ direction.

We update the macroscopic quantities in the Maxwellian by applying a kinetic scheme (see \cite{CP91, Perthame_book} for example) on (\ref{eq:advection_M}). As this equation is solved in an explicit deterministic way, the CFL condition is needed. In the following simulations we find the macroscopic velocity is always less than $5$. So we always take $\Delta t = \Delta x/10$ except in Section \ref{sec:num_conv}. The restriction on the time step due to electric field is $\bigo{1}$, hence satisfied automatically.

The simulations are based on the Landau damping problem. The initial distribution is at equilibrium $f(t=0,x,\bfv) = M(t=0,x,\bfv)$ with the macroscopic quantities given by
\be
\label{num:LD_initial}
\left\{
\ba
&\rho(t=0,x) = 1 + \alpha \sin(x), \\
&\bfu(t=0,x) = 0,\\
& T(t=0,x) = 1,
\ea
\right.
\ee
where $\alpha$ is a constant characterizing how close the initial distribution is to the global equilibrium. Both the linear and nonlinear Landau damping are considered in the following tests.

\subsection{HDP method for VP-BGK system}

In the first test, our aim is to check if our methods can capture the Landau damping properties of the VP-BGK system. In this case the coarse particles are not needed. Thanks to the simple algorithms (\ref{eq:BGK_exp}) or (\ref{eq:BGK_imp}) on the collision, this model allows us to focus on the advection substep.

The reference solution is solved by a PIC-DSMC method.  The advection is solved by a regular PIC method.  The algorithms of Monte Carlo Collision on the homogeneous BGK equation $\partial_t f = Q_{BGK}(f)$ are obtained by applying a forward Euler method
\be
\label{eq:BGK_exp_DSMC}
f^{n+1} = \Delta t \mu M + (1-\Delta t\mu) f^n,
\ee
or a backward Euler method
\be
\label{eq:BGK_imp_DSMC}
f^{n+1} = \frac{1}{1+\Delta t\mu}  M +  \frac{\Delta t \mu}{1+\Delta t\mu}f^n.
\ee
Hence a fraction of $\Delta t\mu$ or $\frac{1}{1+\Delta t\mu}$ particles are removed and replaced by particles sampled from $M$. Again one needs $\Delta t\mu\le1$ in the explicit method (\ref{eq:BGK_exp_DSMC}).

In the following the explicit methods (\ref{eq:BGK_exp}) and (\ref{eq:BGK_exp_DSMC}) are applied. We use the same effective number $\Neff = 5\times10^{-7}$ in both the HDP method and the PIC-DSMC method.

\begin{figure}
\centering
    \includegraphics[width=0.78\textwidth]{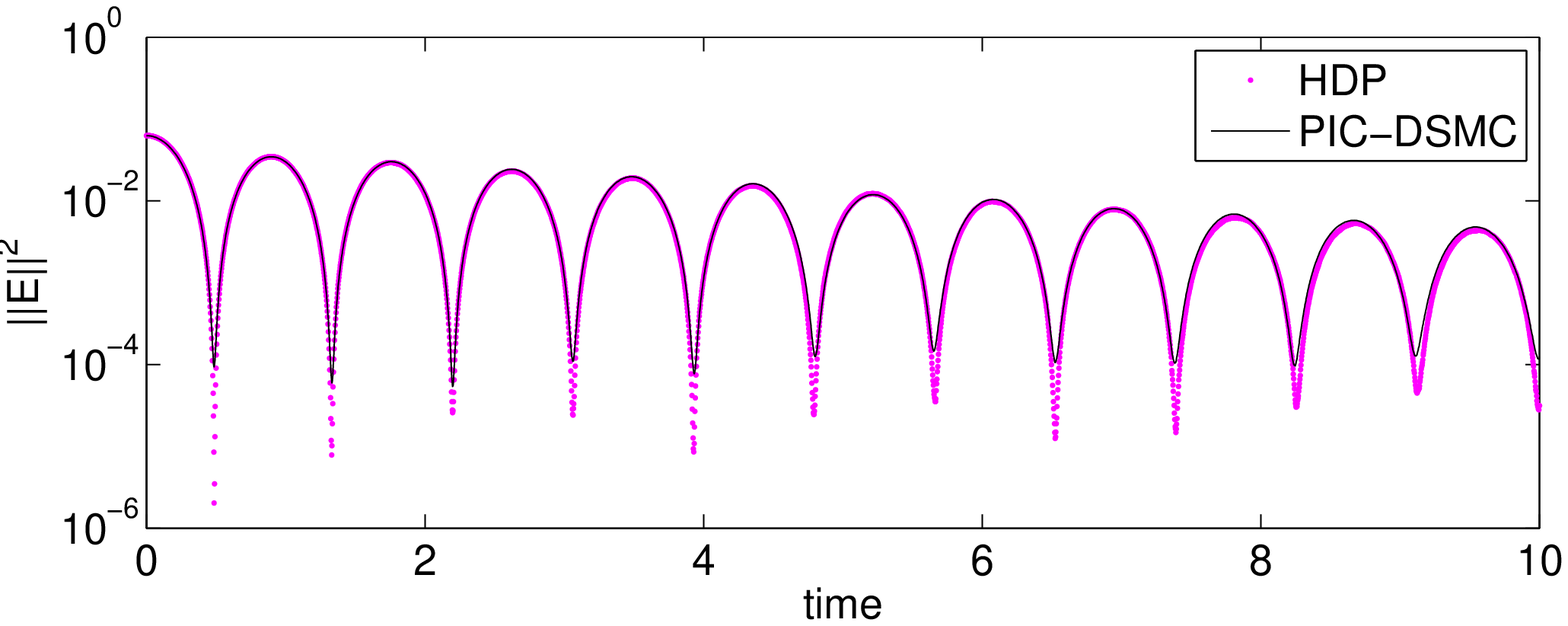}
    \includegraphics[width=0.78\textwidth]{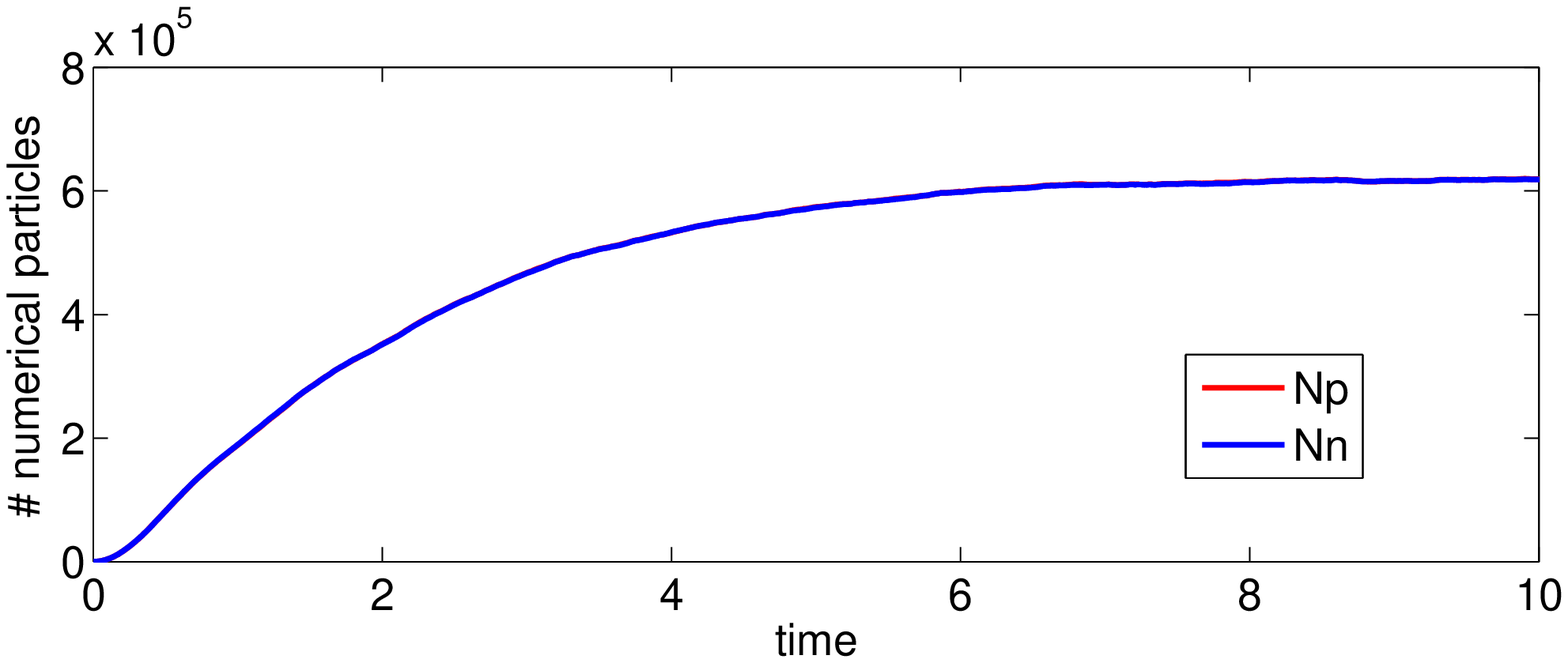}
    \caption{The linear Landau damping problem of the VP-BGK system with collision rate $\mu=1$. Top: the comparison of the decay of energy between the HDP method (red dots) and PIC-DSMC method (black line). Bottom: the growth of the number of positive and negative deviational particles. In the bottom figure, the two lines overlap each other.}
    \label{fig:BGK_lLD}
\end{figure}
\begin{figure}
\centering
    \includegraphics[width=0.78\textwidth]{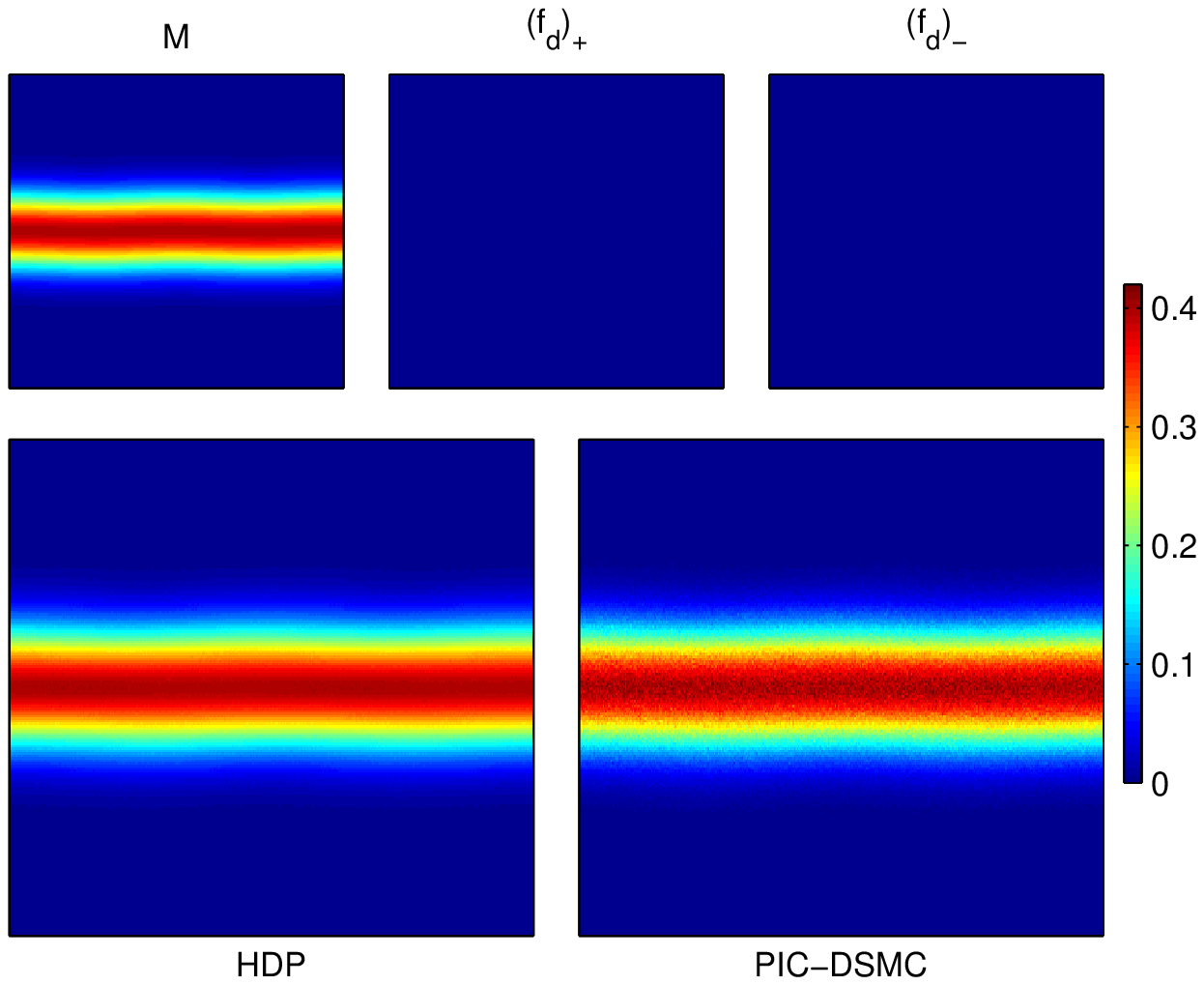}
    \caption{The snapshots of the distribution in the $x-v_1$ phase space at time $t=1.25$ in the linear Landau damping problem of the VP-BGK system with $\mu=1$. The top three figures show the components $M$, $(f_d)_+$ and $(f_d)_-$ in the HDP method. The bottom left figure gives the solution $M+f_d$ obtained by HDP method. The bottom right gives the solution $f$ obtained by PIC-DSMC method.}
    \label{fig:lLD_BGK_snap}
\end{figure}

In Figure \ref{fig:BGK_lLD} and \ref{fig:lLD_BGK_snap} a linear Landau damping problem is studied, with $\alpha = 0.01$ in (\ref{num:LD_initial}). A weak collision frequency $\mu=1$ is used. The top of Figure \ref{fig:BGK_lLD} compares the damping of electric field
\be
\label{eq:elecenergy}
||E||^2 = \Delta x \sum_{k=1}^{N_x} (E_k)^2
\ee
between the solutions obtained by HDP method (red dots) and the PIC-DSMC method (black lines). A good agreement is observed. In the bottom figure, the time evolution of $N_p$ and $N_n$ (i.e. the numbers of positive and negative components of the deviational particles) in HDP method are plotted. The two lines overlap each other, suggesting that the deviation part does not carry any mass. The growth of particles caused by the advection part gets balanced by the collision as time evolves. For comparison, $2.5\times10^{7}$ particles are used in the PIC-DSMC methods. The HDP methods uses less than $3\%$ particles. Figure \ref{fig:lLD_BGK_snap} provides a snapshot of the solutions in the $(x,v_1)$ phase space obtained by HDP method and PIC-DSMC method at time $t=1.25$. The top three figures correspond to $M$, $(f_d)_+$ and $(f_d)_-$ in the HDP solution. The bottom figures correspond to the distribution $f = M + f_d$ in the HDP solution (left) and $f$ in the PIC-DSMC solution (right). In this case the deviation $f_d$ is very small and the solutions obtained in two methods agree well with each other.



\begin{figure}
\centering
    \includegraphics[width=0.78\textwidth]{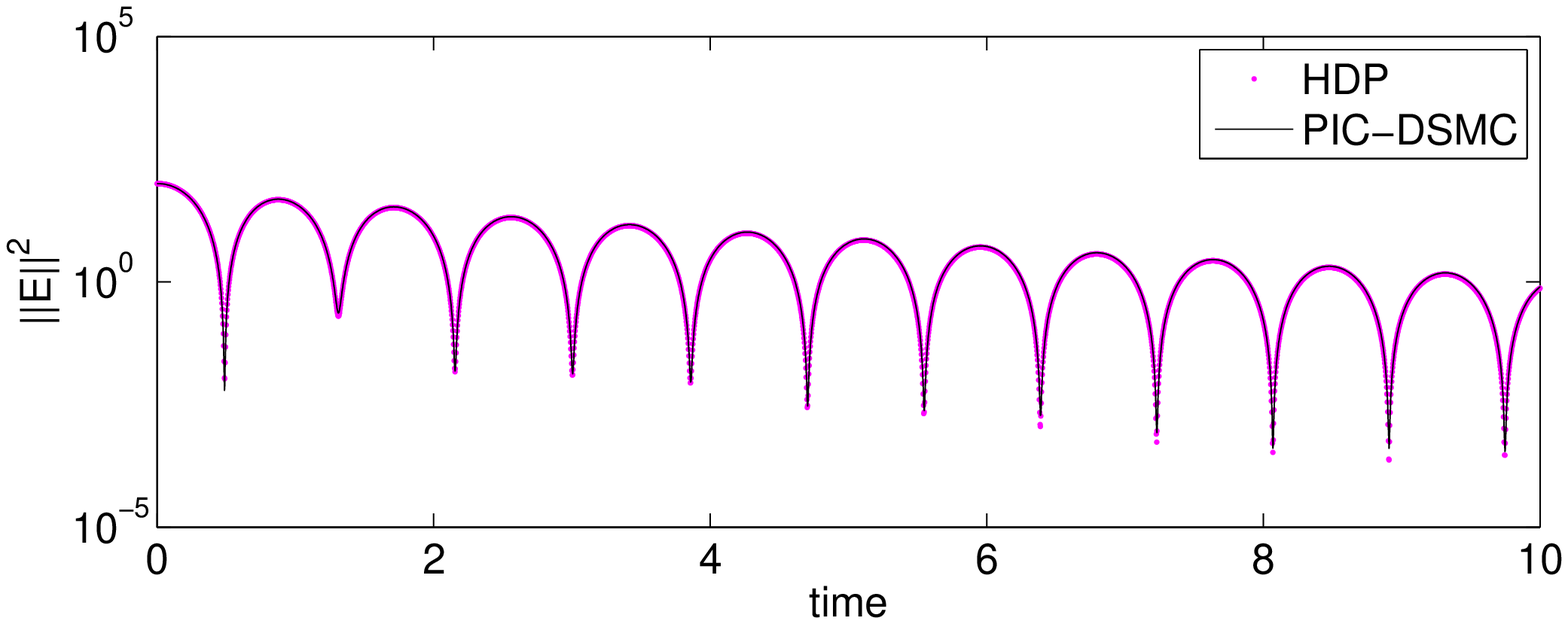}
    \includegraphics[width=0.78\textwidth]{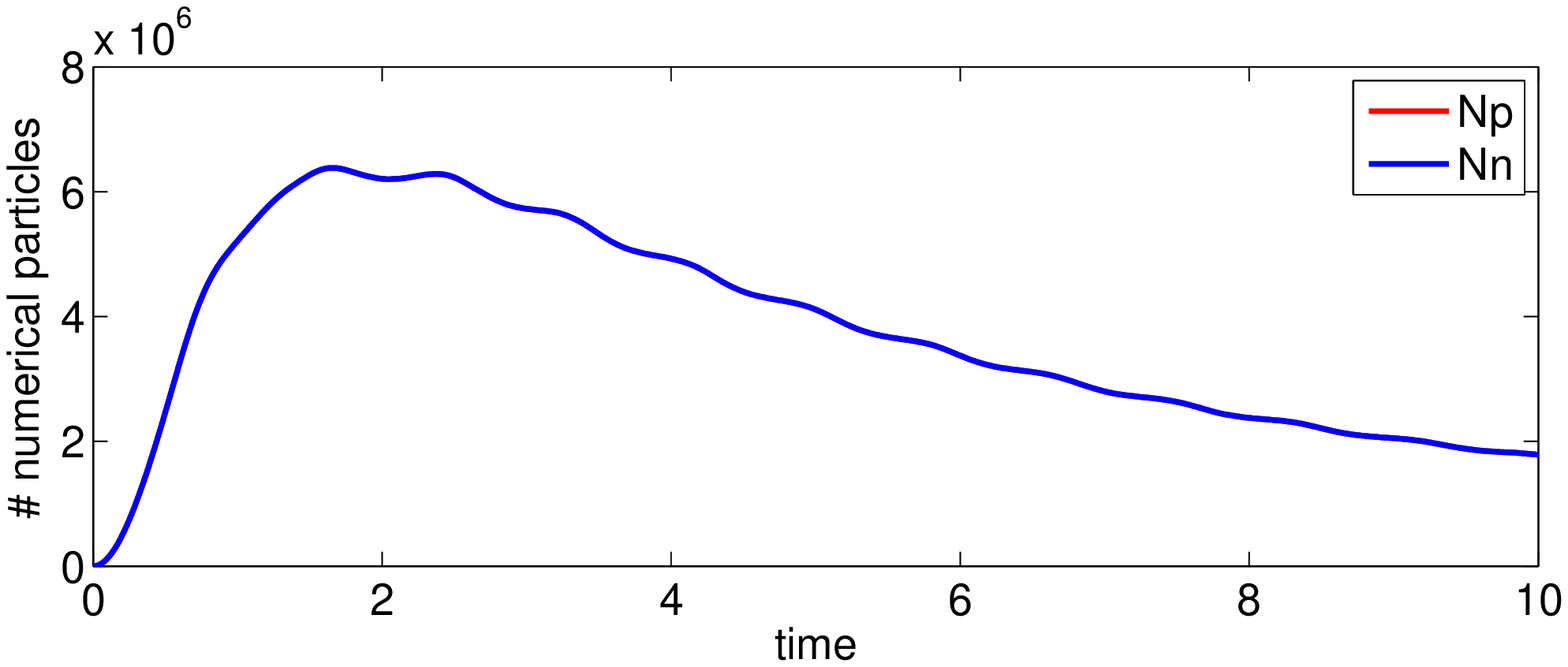}
    \caption{The nonlinear Landau damping problem of the VP-BGK system with collision rate $\mu=1$. Top: the comparison of the decay of energy between the HDP method (red dots) and PIC-DSMC method (black line). Bottom: the growth of the number of positive and negative deviational particles. In the bottom figure, the two lines overlap each other.}
    \label{fig:BGK_nLD}
\end{figure}

\begin{figure}
\centering
    \includegraphics[width=0.78\textwidth]{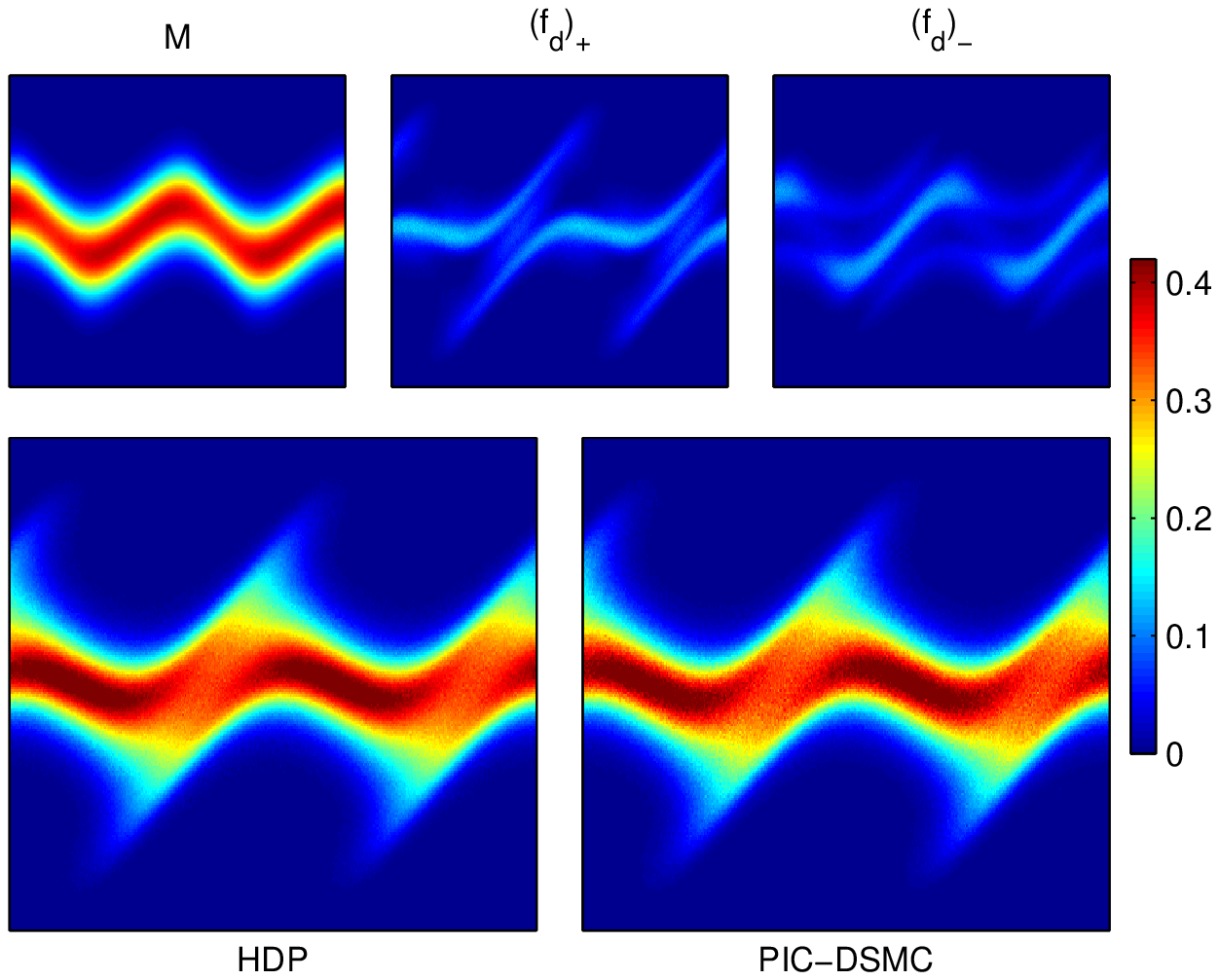}
    \caption{The snapshots of the distribution in the $x-v_1$ phase space at time $t=1.25$ in the nonlinear Landau damping problem of the VP-BGK system  with $\mu=1$. The top three figures show the components $M$, $(f_d)_+$ and $(f_d)_-$ in the HDP method. The bottom left figure gives the solution $M+f_d$ obtained by HDP method. The bottom right gives the solution $f$ obtained by PIC-DSMC method.}
    \label{fig:nLD_BGK_snap}
\end{figure}

In Figure \ref{fig:BGK_nLD} and \ref{fig:nLD_BGK_snap} a nonlinear Landau damping problem is studied, with $\alpha = 0.4$ in (\ref{num:LD_initial}). Again a weak collision frequency $\mu=1$ is used. In this case one can observe the decay of particle numbers as time evolves, which indicates the solution is approaching a global equilibrium which has a uniform density in space. The snapshot in Figure \ref{fig:nLD_BGK_snap} clearly shows the deviation of the solution to the local equilibrium. Again the solutions obtained by the two methods agree well with each other.

\begin{figure}
\centering
    \includegraphics[width=0.78\textwidth]{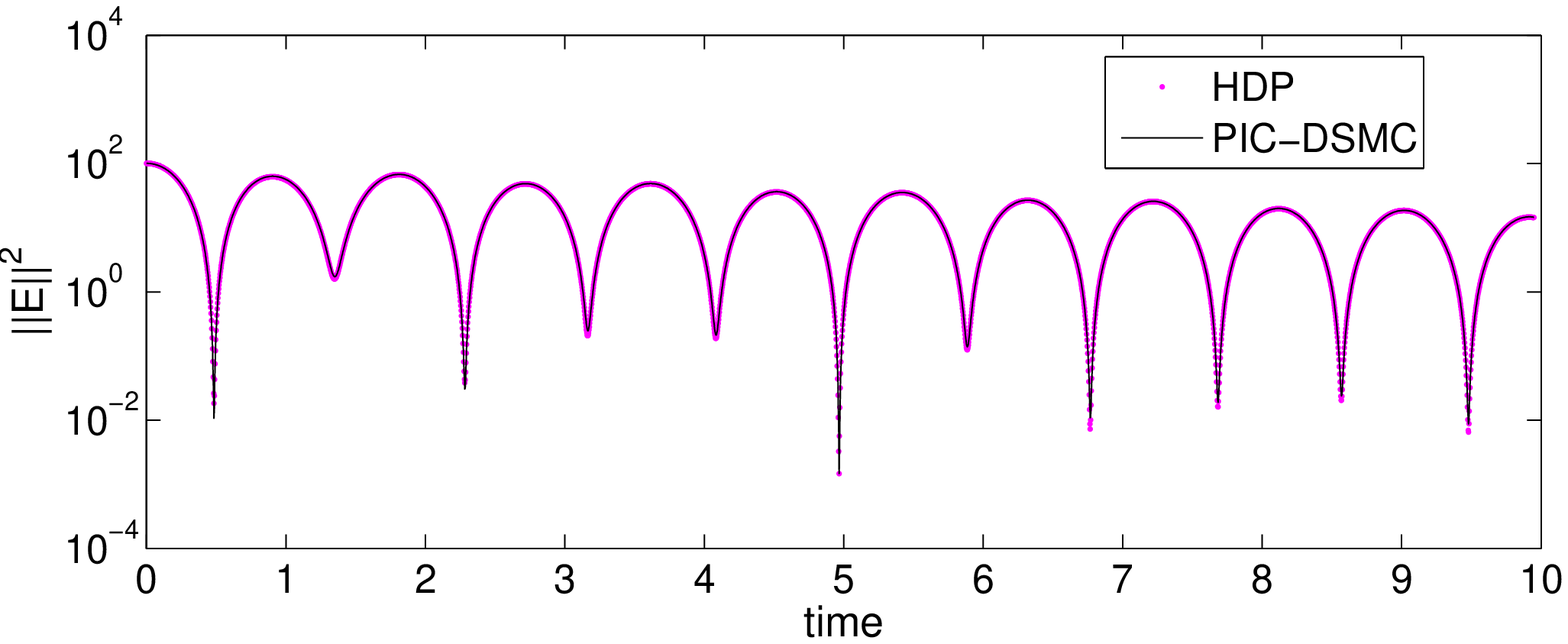}
    \includegraphics[width=0.78\textwidth]{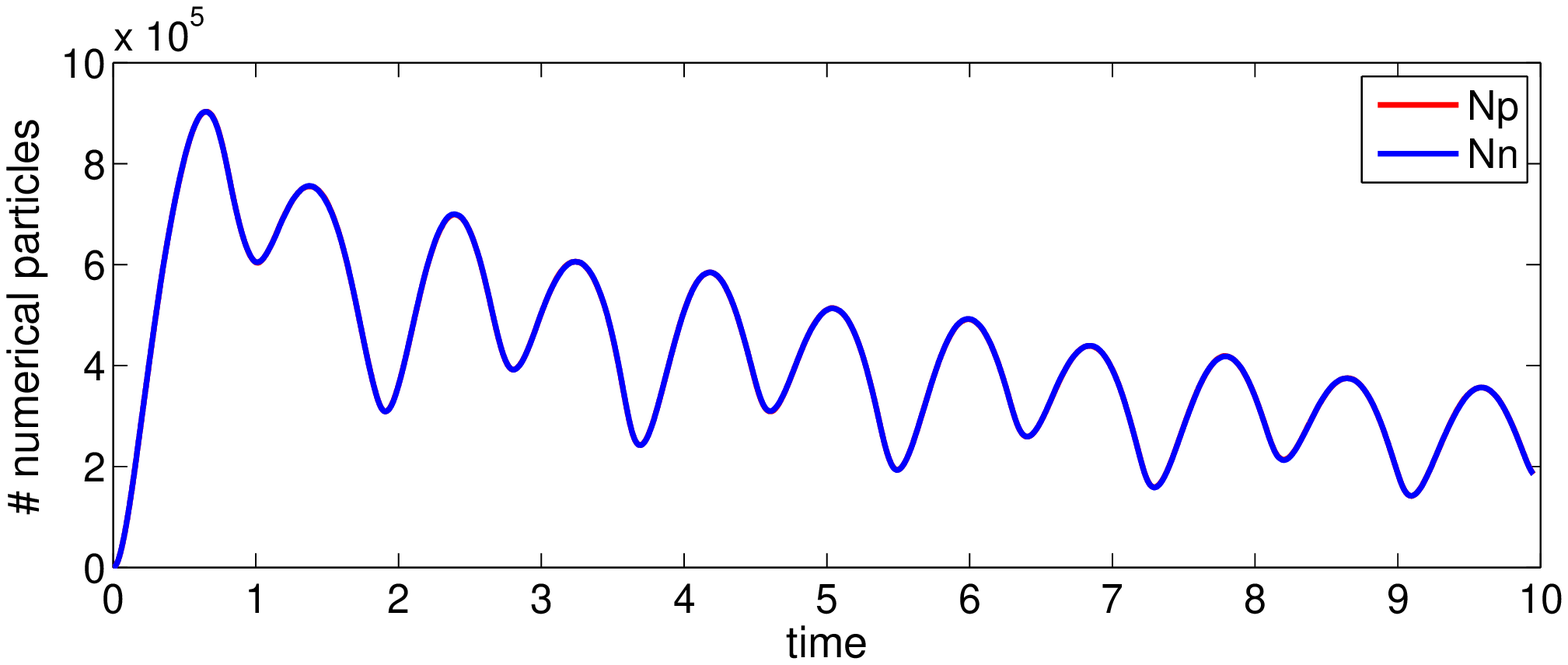}
    \caption{The nonlinear Landau damping problem of the VP-BGK system with a strong collision rate $\mu=10$. Top: the comparison of the decay of energy between the HDP method (red dots) and PIC-DSMC method (black line). Bottom: the growth of the number of positive and negative deviational particles. In the bottom figure, the two lines overlap each other.}
    \label{fig:BGK_nLD_strong}
\end{figure}

\begin{figure}
\centering
    \includegraphics[width=0.78\textwidth]{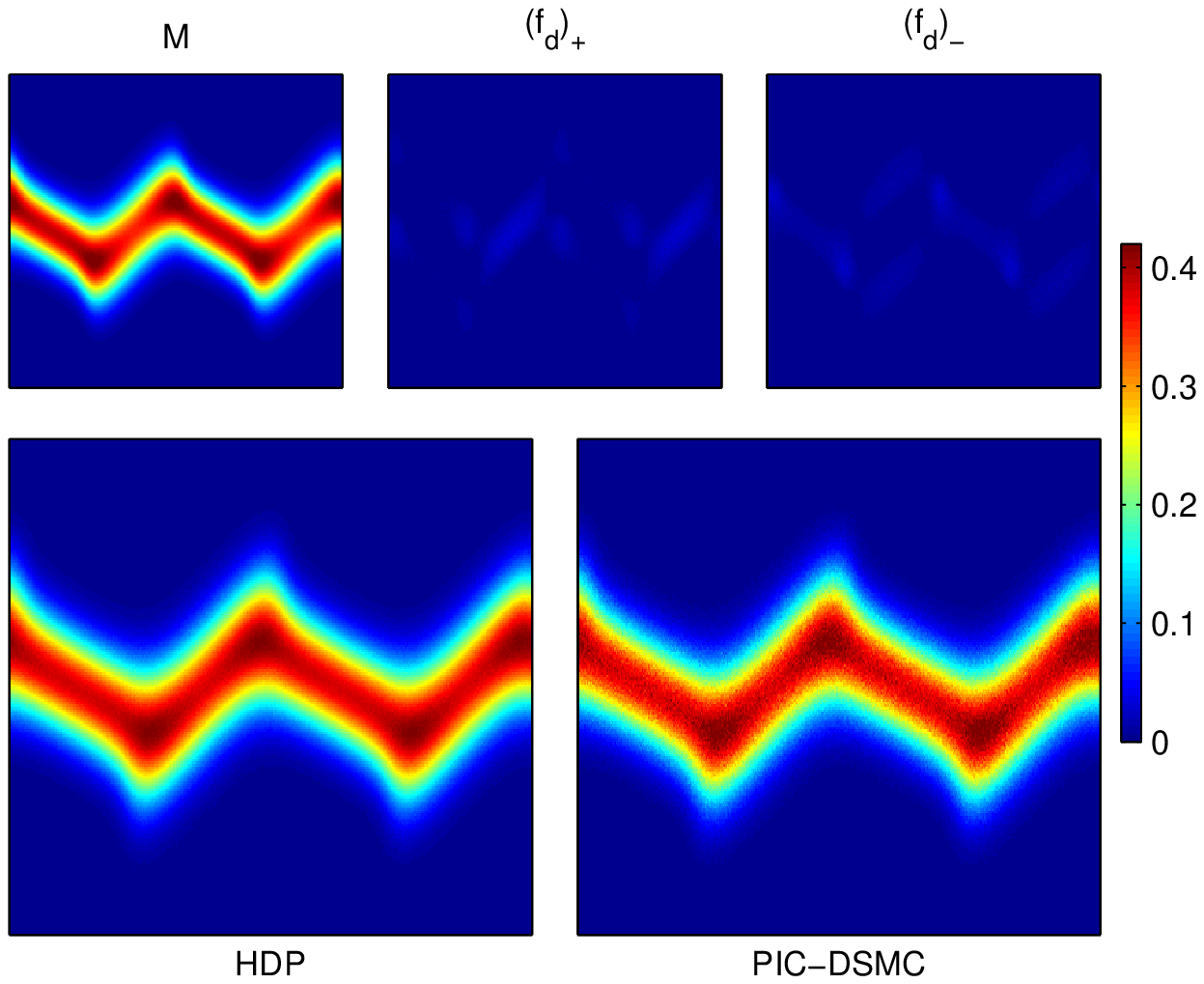}
    \caption{The snapshots of the distribution in the $x-v_1$ phase space at time $t=1.25$ in the nonlinear Landau damping problem of the VP-BGK system  with a strong collision frequency $\mu=10$. The top three figures show the components $M$, $(f_d)_+$ and $(f_d)_-$ in the HDP method. The bottom left figure gives the solution $M+f_d$ obtained by HDP method. The bottom right gives the solution $f$ obtained by PIC-DSMC method.}
    \label{fig:nLD_BGK_strong}
\end{figure}

In Figure \ref{fig:BGK_nLD_strong} and \ref{fig:nLD_BGK_strong} the same nonlinear Landau damping problem is studied with a strong collision frequency $\mu=10$. Now one can observe the oscillation in the particle number, which illustrates the competition between the advection part and the collision part. $N_d$ gradually decays in the long time, which suggests the decay of the solution to the global equilibrium. The deviation part is smaller than the weak collision case and fewer particles are needed.

\subsection{HDP method for VPL system}

Now we apply the HDP method on the VPL system (\ref{eq:VPL}) with the same initial data (\ref{num:LD_initial}). Again the reference solution is solved by a PIC-DSMC method, with a PIC method on the advection, and a DSMC method by Takizuka-Abe \cite{TakizukaAbe} on the Coulomb collisions. In this case we have to include the coarse particles and perform the particle resampling whenever needed.


\subsubsection{Landau Damping}

\begin{figure}
\centering
    \includegraphics[width=0.78\textwidth]{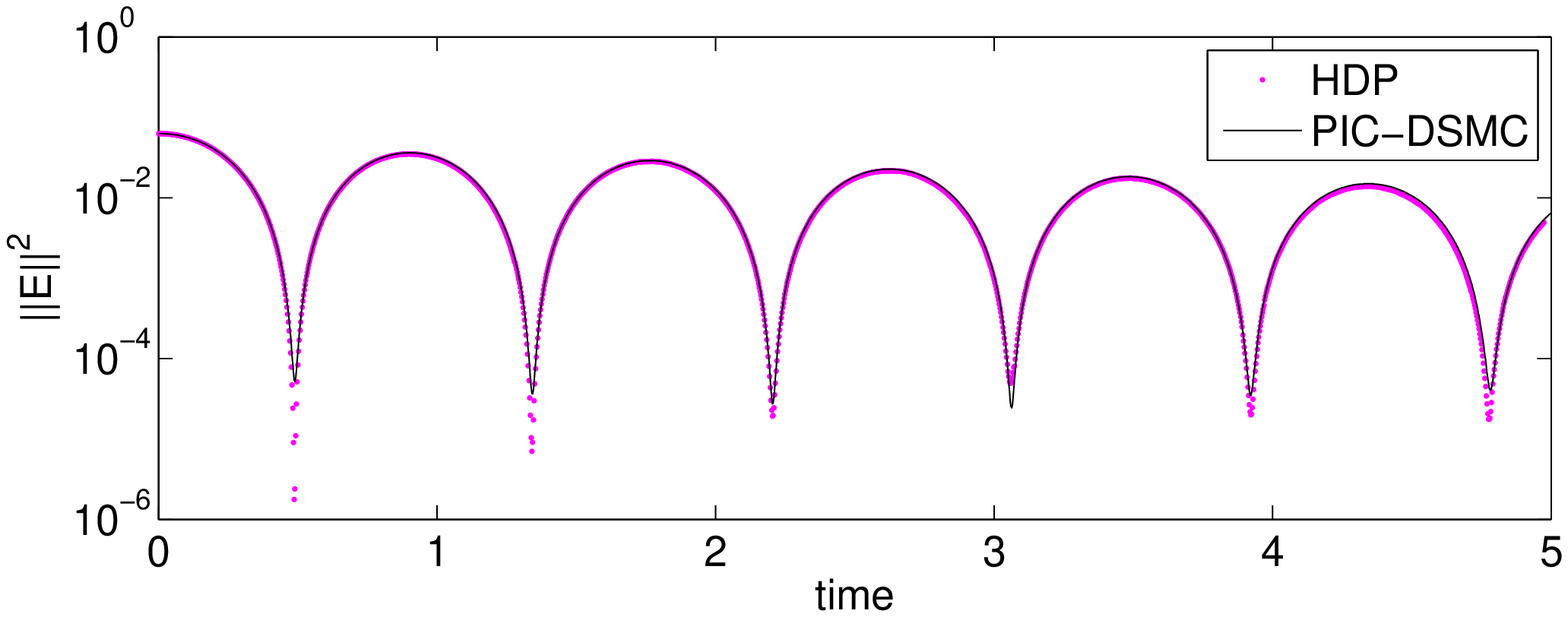} 
    \includegraphics[width=0.78\textwidth]{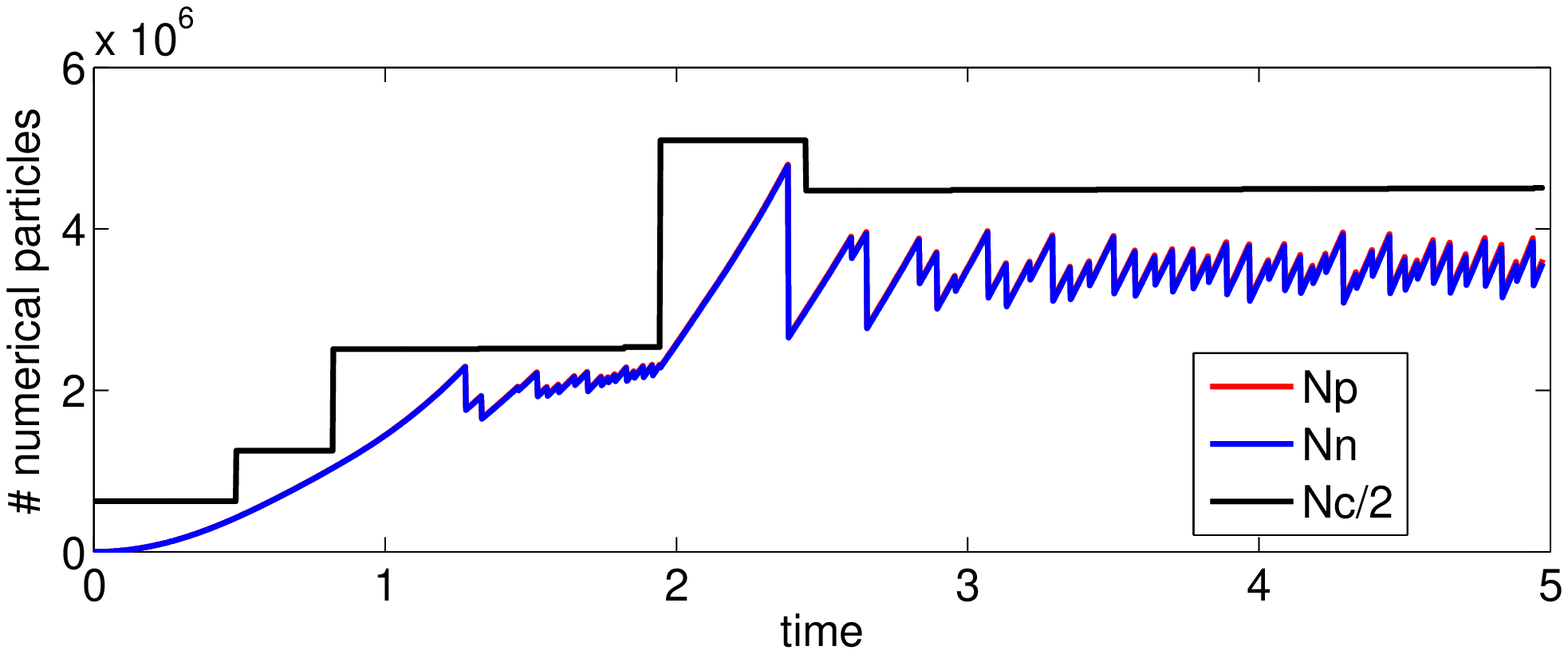}
    \caption{The linear Landau damping problem of the VPL system. Top: the comparison of the decay of energy between the HDP method (red dots) and PIC-DSMC method (black line). Bottom: the growth of the number of positive and negative deviational particles. In the bottom figure, the two lines overlap each other.}
    \label{fig:lLD}
\end{figure}

Figure \ref{fig:lLD} presents the results on the linear Landau damping problem with $\alpha = 0.01$. The top figure shows the good agreement between the solution obtained from the HDP method and the reference solution by a PIC-DSMC method. The bottom figure gives the time evolution of the two components of $N_d$ (the number of deviational particles) and $N_c$ (the number of coarse particles). The zigzag pattern in $N_d$ results from the resampling of coarse particles. In the early stage, one can see the growth of $N_c$ when the resampling of deviational particles fails, which is the case 2 discussed in Section \ref{sec:resample_coarse}.  After the transient time period, $N_d$ and $N_c$ stop from growing. The reduction of $N_c$ from resampling can be observed, which is the case 3 in Section \ref{sec:resample_coarse}.

\begin{figure}
\centering
    \includegraphics[width=.78\textwidth]{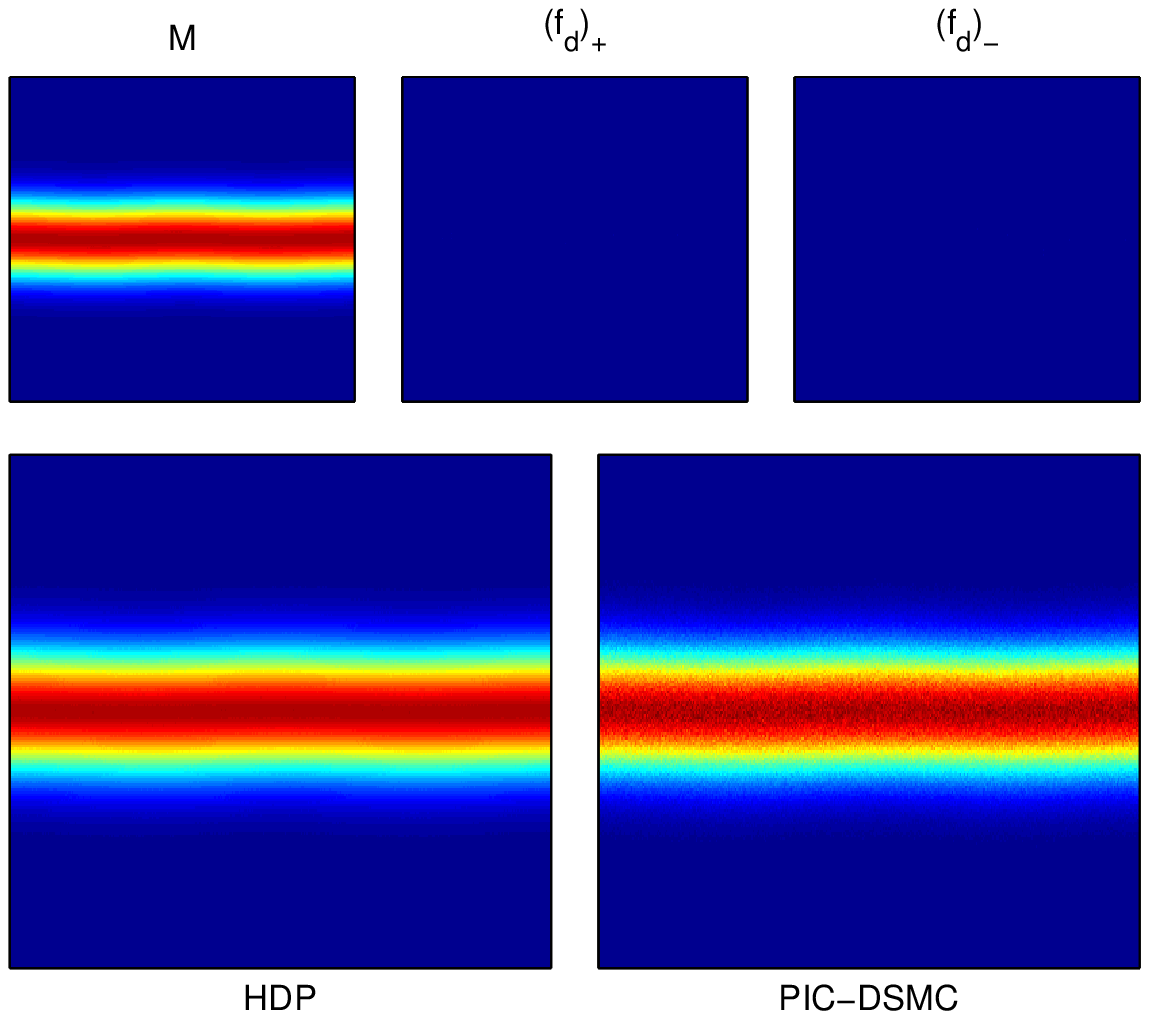}
    \caption{The snapshots of the distribution in the $x-v_1$ phase space at time $t=1.25$ in the linear Landau damping problem of the VPL system. The top three figures show the components $M$, $(f_d)_+$ and $(f_d)_-$ in the HDP method. The bottom left figure gives the solution $M+f_d$ obtained by HDP method. The bottom right gives the solution $f$ obtained by PIC-DSMC method.}
    \label{fig:lLD_snap}
\end{figure}

Figure \ref{fig:lLD_snap} shows the snapshots of the solutions from the HDP method and the PIC-DSMC method at time $t=1.25$. Again, due to the small number of deviational particles, the deviation part $f_d$ is almost invisible. The solutions are close to a global equilibrium.

\begin{figure}
\centering
    \includegraphics[width=0.78\textwidth]{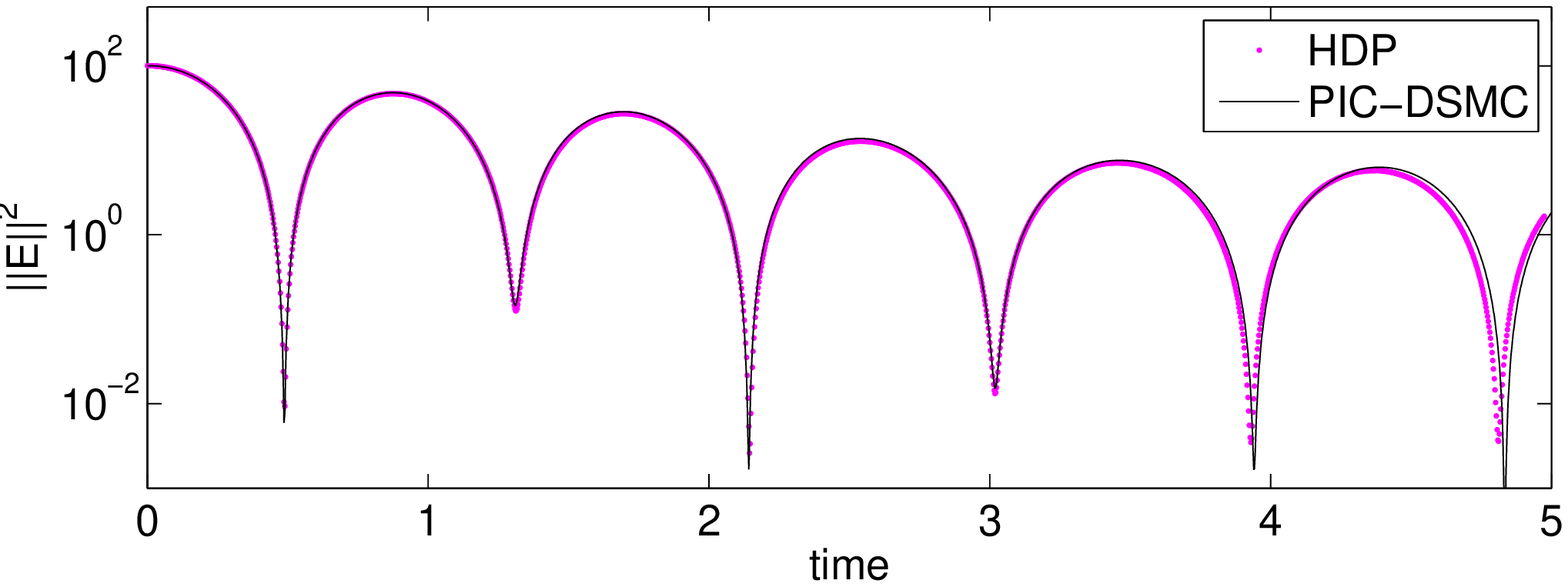} 
    \includegraphics[width=0.78\textwidth]{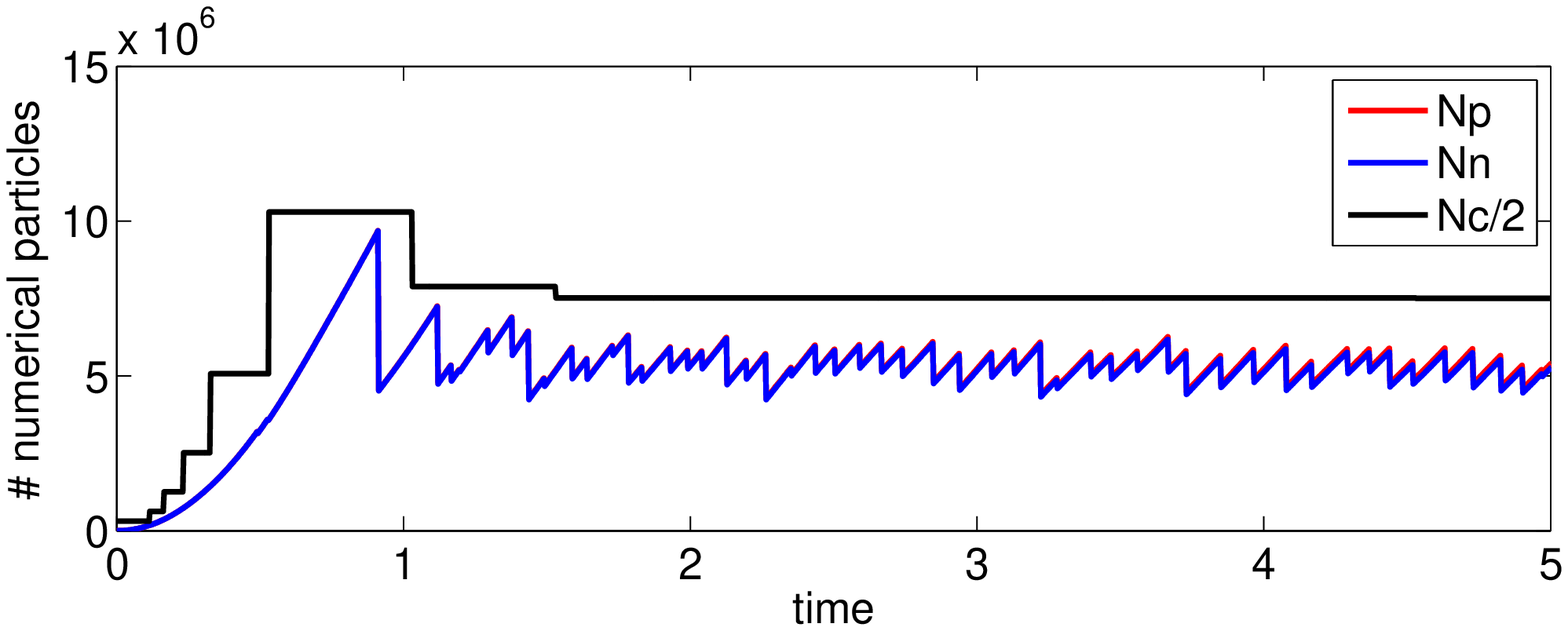} 
    \caption{The nonlinear Landau damping problem of the VPL system. Top: the comparison of the decay of energy between the HDP method (red dots) and PIC-DSMC method (black line). Bottom: the growth of the number of positive and negative deviational particles. In the bottom figure, the two lines overlap each other.}
    \label{fig:nLD}
\end{figure}

Figure \ref{fig:nLD} tests on the nonlinear Landau damping problem with $\alpha = 0.4$. Again, the top figure illustrates the good agreement between the solution from the HDP method and the reference solution. The bottom figure shows the time evolution of $N_d$ and $N_c$. A similar pattern on the growth of the particle numbers is observed as in the linear Landau damping problem.

\begin{figure}
\centering
    \includegraphics[width=.78\textwidth]{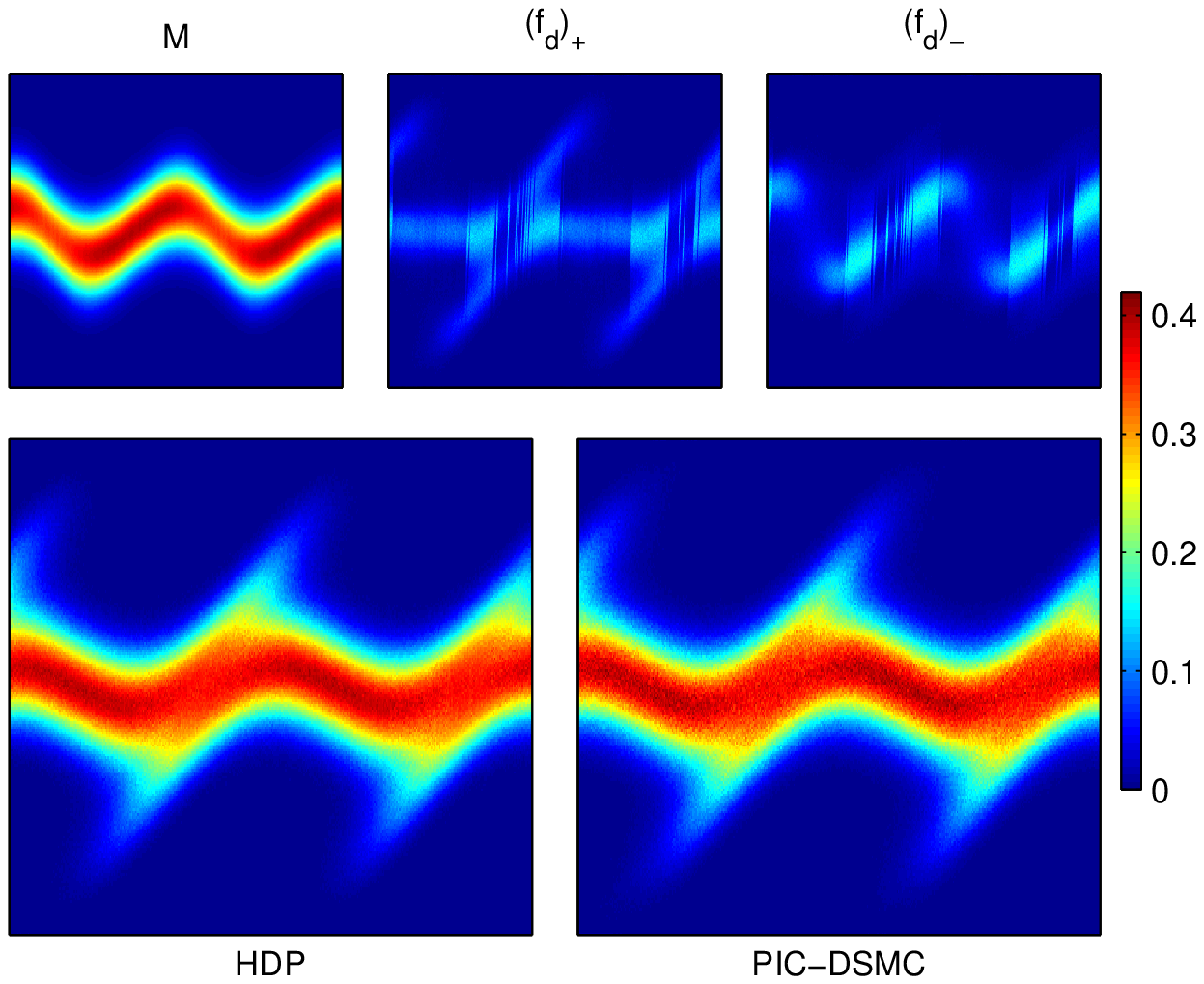}
    \caption{The snapshots of the distribution in the $x-v_1$ phase space at time $t=1.25$ in the nonlinear Landau damping problem of the VPL system. The top three figures show the components $M$, $(f_d)_+$ and $(f_d)_-$ in the HDP method. The bottom left figure gives the solution $M+f_d$ obtained by HDP method. The bottom right gives the solution $f$ obtained by PIC-DSMC method.}
    \label{fig:nLD_snap}
\end{figure}

Figure \ref{fig:nLD_snap} gives the snapshots of the solutions from the HDP method and the PIC-DSMC method at time $t=1.25$. Now the deviation part is significant. The vertical lines in $(f_d)_+$ and $(f_d)_-$ are resulted from the resampling of deviational particles in those cells.

The aim of this simulation is to show that the HDP method is applicable in all regimes. But we would like to mention that the HDP method for nonlinear Landau damping is not more efficient than the PIC-DSMC method since the deviation part $f_d$ is quite large and many deviational particles (and hence coarse particles) are needed.

\subsubsection{Convergence test}
\label{sec:num_conv}
Now we perform a numerical test to check convergence of the HDP method on the VPL system in all regimes, as the effective number $\Neff\to0$ (or equivalently particle number $N_d\to\infty$), and as the time step $\Delta t\to0$.

We consider the relative error in the distribution
\be
\label{def:numeric_error}
\text{error} = \frac{||f-f_\text{ref}||_1}{||f||_1},
\ee
with $||\cdot||_1$ the $L^1$ norm. $f_\text{ref}$ is the reference solution solved by the HDP method with a much smaller time step and effective number.

\begin{figure}
    \includegraphics[width=.48\textwidth]{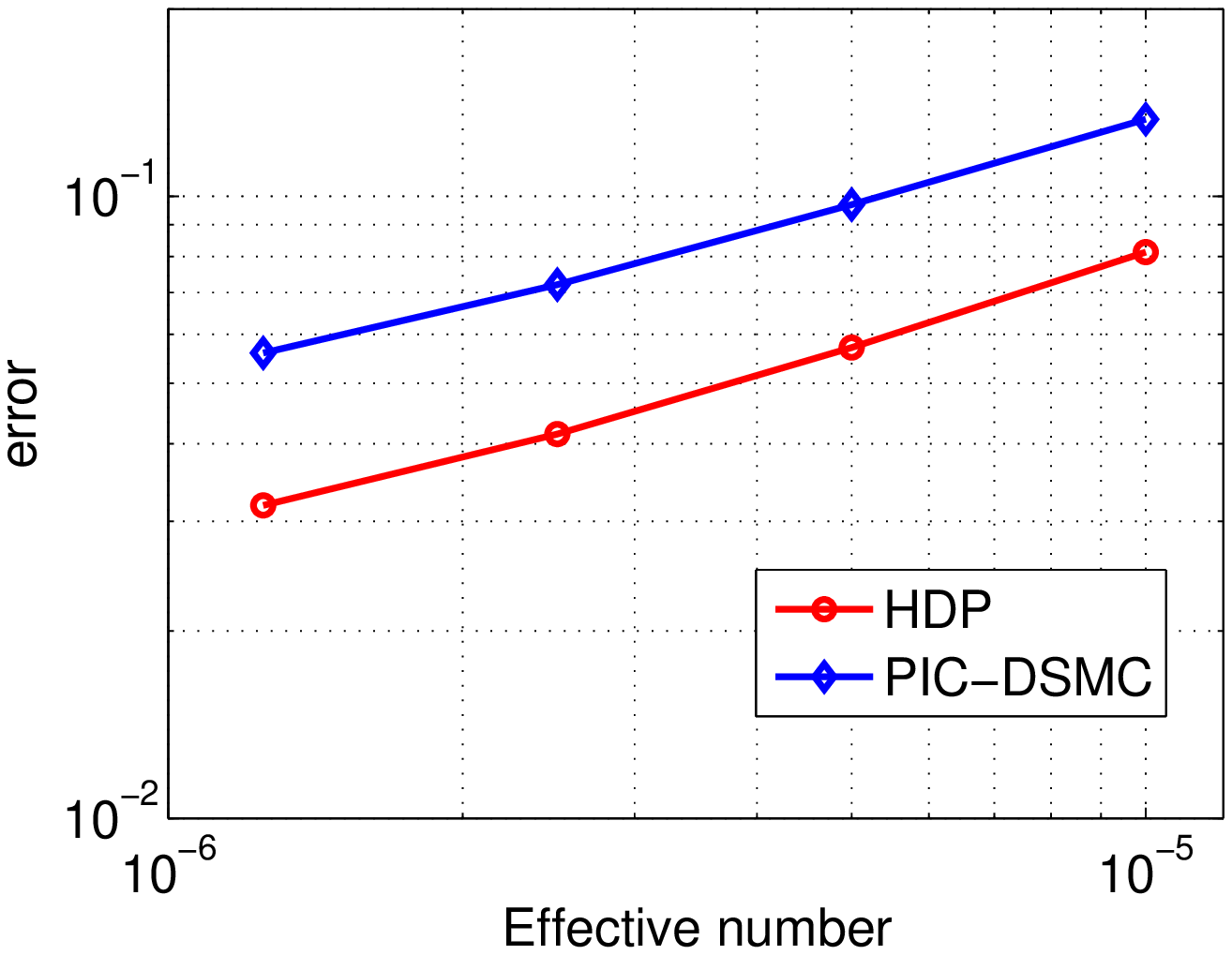}
    \includegraphics[width=.48\textwidth]{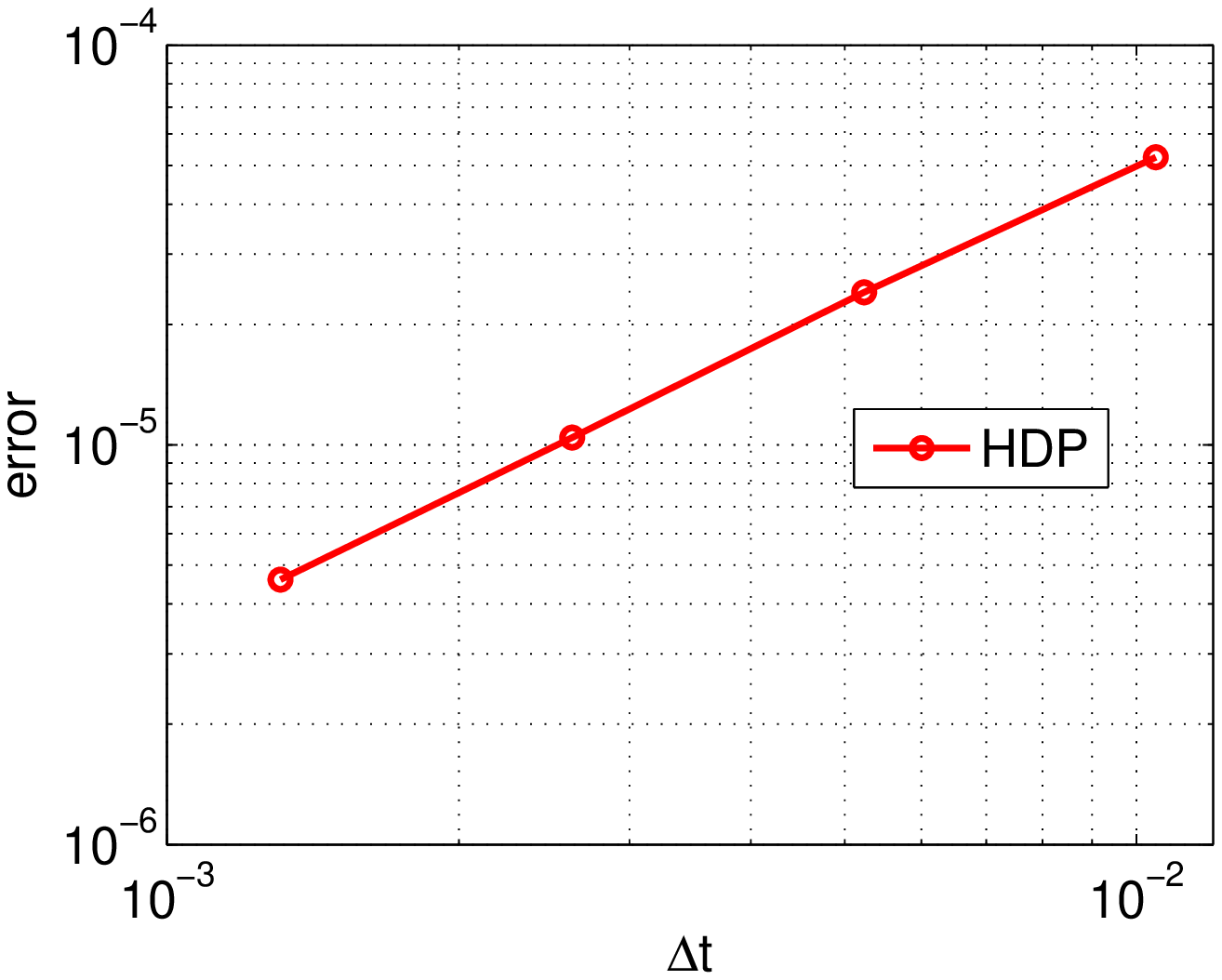}
    \caption{Convergence test  of the HDP method on the VPL system.  Left: the convergence in effective number in nonlinear Landau damping at time $t=0.5$ for the HDP method (red circles) and PIC-DSMC method (blue diamonds). Right: the convergence in $\Delta t$ in linear Landau damping at time $t=0.5$ for the HDP method (red circles).}
    \label{fig:convtest}
\end{figure}

The convergence test as $\Neff\to0$ is performed on the nonlinear Landau damping problem with $\alpha = 0.4$. $\Delta x$ and $\Delta t$ are fixed. We take $\Neff$ ranging from $10^{-5}$ to $1.25\times10^{-6}$. The reference solution is obtained with $\Neff=6.25\times10^{-7}$. The coarse particles use the same effective number. As illustrated in the left of Figure \ref{fig:convtest}, a half order convergence is observed in both the coarse solution, which is solved by a PIC-DSMC method, and the fine solution, which is solved by the HDP method.

The convergence test as $\Delta t\to0$ is performed on the linear Landau damping problem with $\alpha = 0.01$. In this case the effective number is fixed to be $\Neff = 10^{-9}$. As shown in the right of Figure \ref{fig:convtest}, the HDP method has a first order convergence in $\Delta t$. As for the coarse solution solved by the PIC-DSMC method, we find the error in time step $\Delta t$ is completely dominated by the statistical error due to finite number of samples. To observe the first order time step error, the size of samples becomes too large to be simulated. The comparison between the two figures in Figure \ref{fig:convtest} also suggests that the deterministic error (due to time step $\Delta t$) is much smaller than the statistical error (due to particle number $N$).

\subsubsection{Efficiency test}
\label{sec:num_effi}

Finally we perform an efficiency test on the HDP methods for the VPL system. The HDP methods are designed to accelerate the simulation of PIC-DSMC methods in the near fluid regime. To test this, we consider the Landau damping problem with different $\alpha$'s. For each $\alpha$, several different effective numbers are used in computing the error (\ref{def:numeric_error}). $\Delta x$ and $\Delta t$ are fixed in this test.

\begin{figure}
\centering
    \includegraphics[width=.95\textwidth]{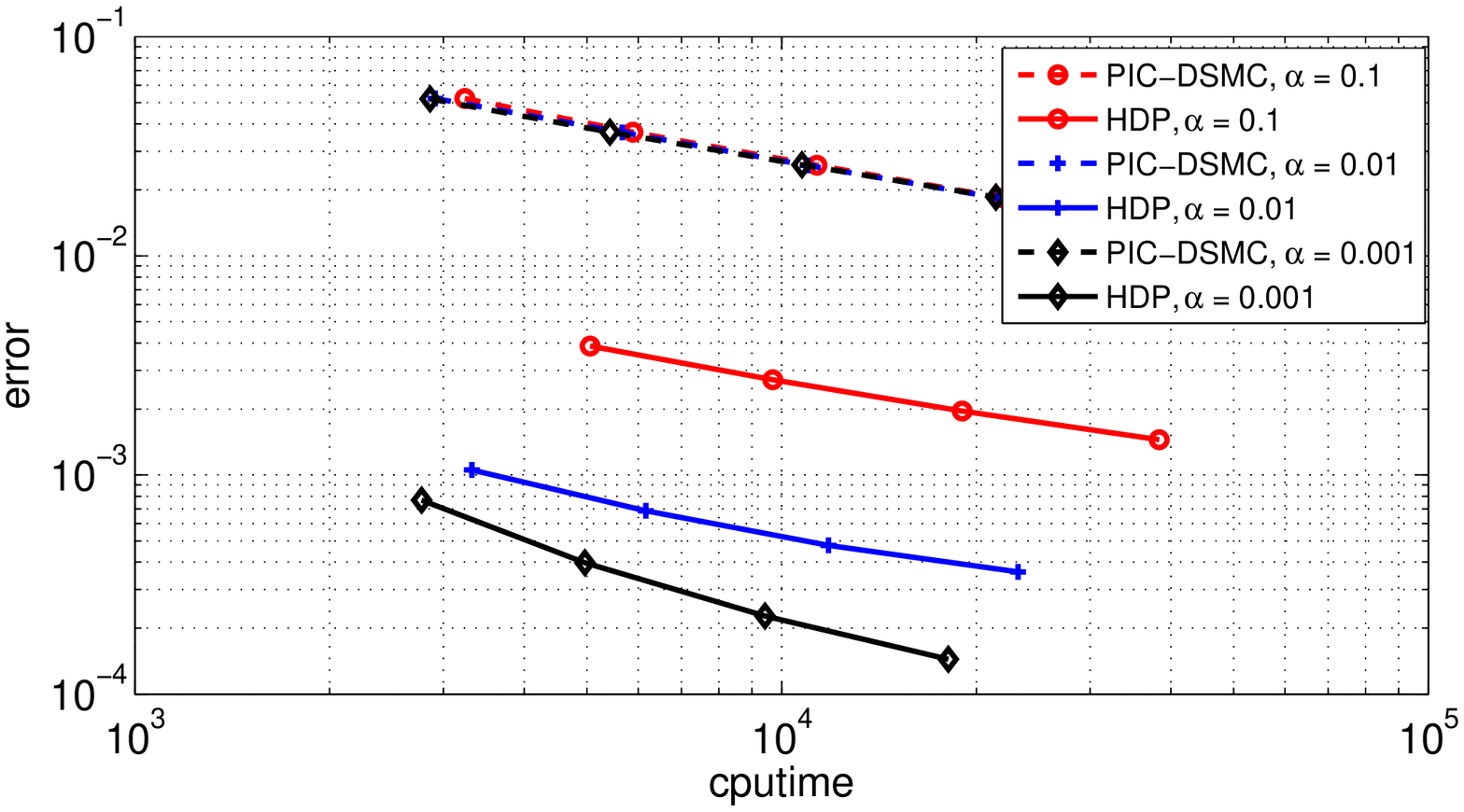}
    \caption{The efficiency test of the HDP method on the VPL system. The initial data (\ref{num:LD_initial}) is used with different $\alpha$. The dash lines represent the PIC-DSMC methods while the solid lines correspond to the HDP methods. }
    \label{fig:effi}
\end{figure}

In Figure \ref{fig:effi} we compare the error against the cpu time in simulation. Solutions by PIC-DSMC methods (dashed lines) and HDP methods (solid lines) with $\alpha = 0.1, 0.01$ and $0.001$ are computed. The three dashed lines almost overlap, suggesting that the efficiency of PIC-DSMC methods is not sensitive on whether the solution is close to the equilibrium. The solid lines give much smaller errors with the same cpu costs, illustrating the high efficiency of the HDP methods. Furthermore, Figure \ref{fig:effi} shows the HDP methods become more efficient when the solution is closer to the equilibrium. This is due to the fact that less and less deviational and hence coarse particles are needed. The HDP methods effectively accelerate the PIC-DSMC methods in the near fluid regime.

\section{Conclusion and Future work}

We have proposed a Hybrid method with Deviational Particles (HDP) to solve the inhomogeneous Vlasov-Poisson-Landau system. This method combines several numerical techniques. The Maxwellian part is solved by a grid based fluid solver. To evolve the deviation part, which is approximated by (positive and negative) deviational particles, we apply the negative particle based Monte Carlo method proposed in \cite{YC15} on the collision part and a Macro-Micro projection in \cite{BLM08} on the advection part. The coarse particles are simulated by a traditional PIC-DSMC method. Moreover, the particle resampling techniques are further investigated and improved.

A next step is the extension of this method to mixed regimes in which the collision frequency could range over several orders of magnitude in space. There could exist both fluid regime (where $f\approx M$) and kinetic regime (where $f$ is far away from $M$). Ideally one uses only a small number of deviational particles near the fluid regime. However in current HDP method the coarse particles have the same weights in all spatial cells. The global effective number for the coarse particles is determined by the kinetic regime where a lot of coarse particles are needed. Hence in the fluid regime, in which the same effective number is used, one also needs to evolve many coarse particles. An efficient method to transport particles with spatial variant effective numbers is to be studied in the future work.

Furthermore, as the coarse solution $f_c$ in the HDP method for VPL system is solved by a traditional PIC-DSMC method, $f_c$ itself is a good approximation of the solution if it is far from Maxwellian. Hence in the kinetic regime, one can discard the deviational solution $(M+f_d)$, and simply use $f_c$ as the fine solution. At the interface of the two regimes, a resampling process can convert between the deviational particles and the coarse particles. This will be further investigated.

\section{Acknowledgement}

The author would like to  gratefully and sincerely thank Professor Russel Caflisch in University of California, Los Angeles for his support, guidance and fruitful discussion during this work.

\bibliography{montecarlo}

\bibliographystyle{plain}

\appendix

\section{Compute the source term in (\ref{eq:advection_fd})}
\label{append:advect_source}

In this appendix we show how to efficiently sample deviational particles from the source term in the right hand side of (\ref{eq:advection_fd}).

Suppose we have the moments of the Maxwellian at time $t^n$ and the electric filed $\bfE$ has been solved from the Poisson equation (\ref{eq:VPL_E}). Then the term $\mathcal T M$  is just the Maxwellian $M$ multiplied by a cubic polynomial in $\bfv$,
\[
\ba
\frac{\mathcal T M}{M} &= \frac{1}{M} (\bfv\cdot\nabla_\bfx M - \bfE\cdot \nabla_\bfv M) \\
& = \bfv^T \left(\frac{\nabla_\bfx \rho}{\rho} + \frac{\nabla_\bfx T}{2 T} \left(\frac{|\bfv-\bfu|^2}{T} - d\right) + \nabla_\bfx^T \bfu \frac{\bfv-\bfu}{T}\right) + \bfE^T \cdot \frac{\bfv-\bfu}{T},
\ea
\]
with the matrix
\[\nabla^T \bfu  = [\nabla_\bfx u_1, \nabla_\bfx u_2, \nabla_\bfx u_3].\]
Here $\rho$, $\bfu$ and $T$ are the density, macroscopic velocity and temperature of the Maxwellian part.

To compute $\Pi_M \mathcal T M$, from the definition (\ref{eq:projection}), one needs to evaluate the the following coefficients
\[
\ba
&\left< \mathcal T M\right> = \nabla_\bfx \cdot (\rho \bfu), \\
&\left< (\bfv-\bfu)\mathcal T M\right> = \rho (\bfu\cdot \nabla_\bfx) \bfu + \nabla_\bfx (\rho T) + \rho \bfE, \\
&\left< \left(\frac{|\bfv-\bfu|^2}{T} - d\right) \mathcal T M\right> = \frac{\rho d}{T} \bfu\cdot \nabla_\bfx T + 2\rho \nabla_\bfx \cdot \bfu. \\
\ea
\]

Similarly, for the term $\Pi_M (\mathcal T f_d)$, one needs to evaluate the following moments for $\phi = 1, \bfv,|\bfv|^2$,
\[\left< \phi\mathcal T f_d \right> = \nabla_x\cdot\left< \bfv f_d\phi\right>,\]
where we have used $\left<\phi\bfE\cdot \nabla_v f_d\right>=0$. Note that this has been computed in (\ref{eq:vfphi}) to update $M$ via (\ref{eq:advection_M}) and can be reused.

Since the projection $\Pi_M \psi$ is just the Maxwellian $M$ multiplied by a quadratic function in $\bfv$, the whole source term in (\ref{eq:advection_fd}) is in the form of $P_3(\bfv) M(\bfv)$, with $P_3(\bfv)$ some cubic function in $\bfv$. One can easily find another Maxwellian with slightly larger temperature as an upper bound. Then deviational particles can be sampled efficiently by a regular rejection-acceptance method.

\section{Acceleration of Particle Resampling}
\label{append:resample_acc}

In this section we describe an approximation of the resampling formula (\ref{eq:resample_formula}) which greatly accelerates the resampling process. We also test on the efficiency numerically.

\subsection{Acceleration Technique}
We use the Fourier basis $\phi_\bfk(\bfv) = \exp(-i\bfk\cdot\bfv)$ in the expansion (\ref{eq:resample_formula}). Note that the evaluation of each Fourier coefficient
\[\hat f_\bfk = \Neff \left( \sum_{\bfv_{p}} \phi_\bfk (\bfv_{p}) - \sum_{\bfv_{n}} \phi_\bfk (\bfv_{n}) \right) \]
needs to compute $\bigo{N_d}$ exponents. This leads to the expensive cost of $\bigo{N_dK^3}$  to evaluate all Fourier coefficients in the 3-dimensional velocity space. In the following we give a much cheaper approximation of $\hat f_\bfk$.

After shifting and rescaling, we can assume all the deviational particles are located in the cube $[0, 2\pi]^3$ in the velocity space. Place a uniform grid on $[0, 2\pi]^3$ with $\tilde K (\ge K)$ points in each direction. For each positive deviational particle $\bfv_p$, write it as
 \[\bfv_p = \bfv_{\xi(p)} + \delta_p,\]
 where $\bfv_{\xi(p)}$ is the nearest grid point of $\bfv_p$. Let $I_\bfj = \{p : \xi(p) = \bfj\}$, which is the set of all particles near the $\bfj$th grid. Define the following functions on the grids
 \[
 \ba
 f^{(0)}(\bfv_\bfj) &= \Neff \sum_{p\in I_\bfj} 1, \\
 f^{(1)}(\bfv_\bfj) &= \Neff \sum_{p\in I_\bfj} \delta_p, \\
 f^{(2)}(\bfv_\bfj) &= \Neff \sum_{p\in I_\bfj} \delta_p \delta_p^T. \\
 \ea
 \]
 Note that the vector $ f^{(1)}$ has $3$ components and the matrix $f^{(2)}$ has $3^2$ components. Then the first term in the coefficient of (\ref{eq:resample_formula}) is
\be
\label{eq:resample_formula_approx}
\ba
 &\Neff\sum_{\bfv_{p}} \exp(-i \bfk \cdot \bfv_{p}) \\
 & = \Neff\sum_{\bfv_{p}} \exp(-i \bfk \cdot \bfv_{\xi(p)})\exp(-i \bfk \cdot \delta_{p}) \\
 & =  \Neff\sum_{\bfv_{p}} \exp(-i \bfk \cdot \bfv_{\xi(p)}) \left( 1 - i \bfk \cdot \delta_p  - \frac{1}{2}(\bfk \cdot \delta_p)^2 \right) + \bigo{(|\bfk|\Delta v/2)^3} \\
 & =  \Neff\sum_{\bfj} \exp(-i \bfk \cdot \bfv_{\bfj}) \sum_{p\in I_\bfj} \left(1 - i \bfk \cdot \delta_p  - \frac{1}{2}(\bfk \cdot \delta_p)^2 \right)+ \bigo{(|\bfk|\Delta v/2)^3} \\
 &= \sum_{\bfj} \exp(-i \bfk \cdot \bfv_{\bfj}) \left( f^{(0)}(\bfv_\bfj) -i \bfk\cdot f^{(1)}(\bfv_\bfj) -\frac{1}{2} \bfk^T f^{(2)}(\bfv_\bfj) \bfk \right) + \bigo{(|\bfk|\Delta v/2)^3} \\
 & = \hat f^{(0)}_\bfk - i \bfk \hat f^{(1)}_\bfk -\frac{1}{2} \bfk^T \hat f^{(2)}_\bfk \bfk+ \bigo{(|\bfk|\Delta v/2)^3}.
  \ea
 \ee

 Therefore one can first recover the functions $f^{(0)}$, $f^{(1)}$ and $f^{(2)}$ on the grids, which takes a number of $\bigo{N_d}$ computations of polynomials. Then apply a Fast Fourier Transform to each function and multiply by $\bfk$ and $\bfk\bfk^T$, which costs $\bigo{\tilde K^3 \log \tilde K}$. Only a cost of $\bigo{N_d + \tilde K^3 \log \tilde K}$ is needed to get an approximation. In practice we take $\tilde K = 2K$. With $K=30$, the omitted high order term is $\bigo{(\tilde K)^{-3}}\approx 10^{-6}$, when evaluating the low modes.

 The second term in the coefficient of (\ref{eq:resample_formula}) can be approximated in a similar way on the negative deviational particles. Once the explicit formula (\ref{eq:resample_formula}) is recovered, one can use a similar approximation to evaluate the function $f(\bfv)$ for each sampled particle with velocity $\bfv$.

\subsection{Numerical Test}

\begin{figure}
\centering
    \includegraphics[width=.95\textwidth]{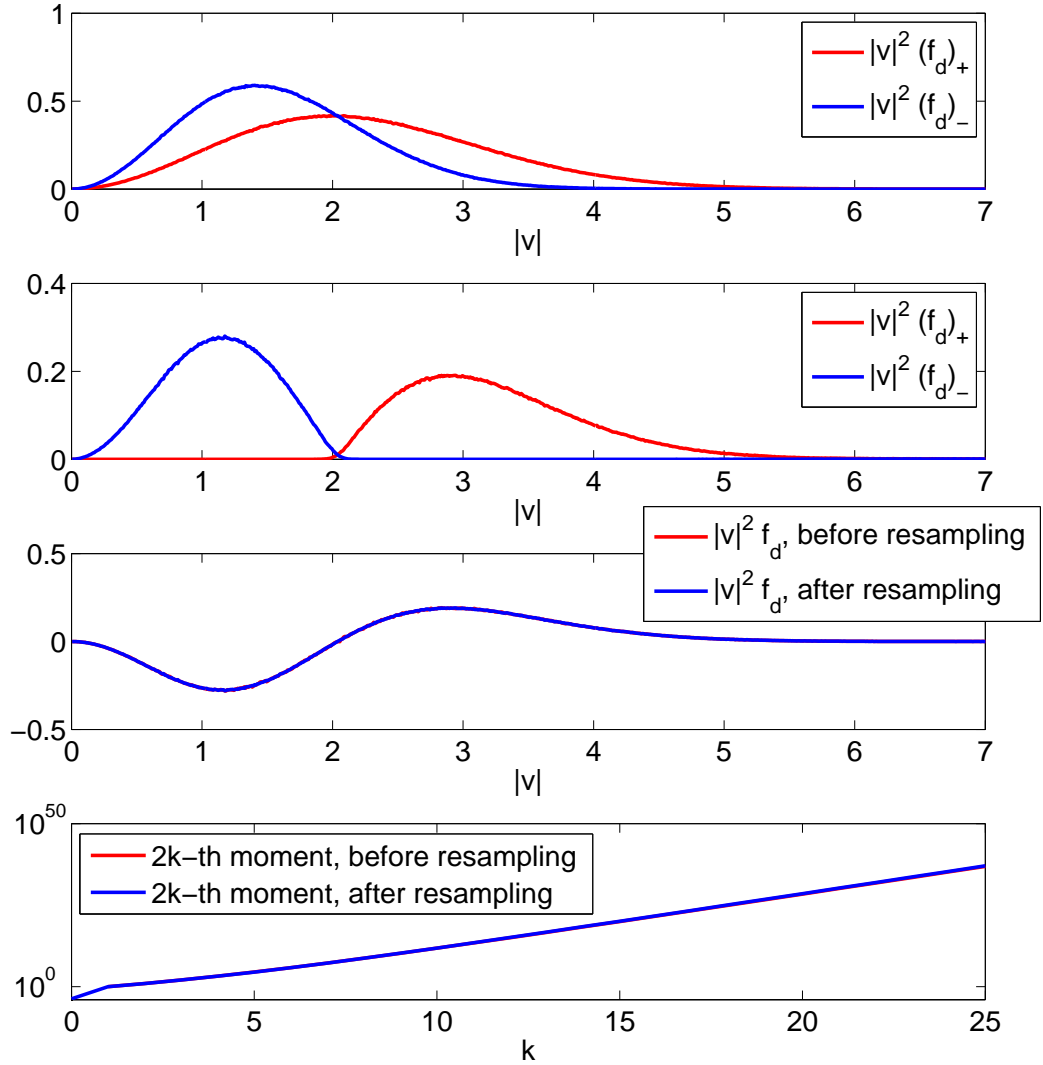}
    \caption{The resampling with the approximation (\ref{eq:resample_formula_approx}). From top to bottom: the distributions $|\bfv|^2 (f_d)_+$ and $|\bfv|^2 (f_d)_-$ before resampling; the distributions $|\bfv|^2 (f_d)_+$ and $|\bfv|^2 (f_d)_-$ after resampling; the distribution $|\bfv|^2 f_d$ before and after resampling; the moments $\int (v_1)^{2k} f_d(\bfv) \ud \bfv$ before and after resampling.}
    \label{fig:resample_fast}
\end{figure}

We perform a numerical test on the resampling of deviational particles in a single cell. The positive and negative deviational particles are sampled from two Maxwellians in 3D velocity space with
\[\rho_p=1, \bfu_p =0, T_p = 2, \rho_n=1, \bfu_n =0, T_n = 1.\]
$10^7$ particles are sampled from each distribution. $K=30$ modes are used. Figure \ref{fig:resample_fast} shows the distribution $|\bfv|^2f_d(\bfv)$ in the energy space and the $2k$-th moments $\int (v_1)^{2k} f_d(\bfv) \ud \bfv$ before and after resampling with the approximation (\ref{eq:resample_formula_approx}). A good agreement is observed. This single resampling process takes about 20 seconds by a serial C++ code on a single core 2.0 GHz CPU processor. As for the resampling with the full expansion (\ref{eq:resample_formula}), we find it cannot be finished in $24$ hours on the same processor if $10^7$ deviational particles are used. For comparison, with $10^4$ deviational particles, the resampling with the full expansion takes about 5 minutes.

\end{document}